\documentclass[10pt]{amsart}
\usepackage[cp1251]{inputenc}
\usepackage[english]{babel}
\usepackage{amsmath}
\usepackage{amssymb}
\usepackage{amsfonts}

\PassOptionsToPackage{linktocpage}{hyperref}

\begin{document}
\renewcommand{\refname}{References}

\thispagestyle{empty}

\title[Expansion of Iterated Stochastic Integrals]
{Expansion of Iterated Stochastic Integrals 
with Respect to Martingale Poisson Measures and 
with Respect to Martingales
Based on Generalized Multiple Fourier Series}
\author[D.F. Kuznetsov]{Dmitriy F. Kuznetsov}
\address{Dmitriy Feliksovich Kuznetsov
\newline\hphantom{iii} Peter the Great Saint-Petersburg Polytechnic University,
\newline\hphantom{iii} Polytechnicheskaya ul., 29,
\newline\hphantom{iii} 195251, Saint-Petersburg, Russia}%
\email{sde\_kuznetsov@inbox.ru}
\thanks{\sc Mathematics Subject Classification: 60H05, 60H10, 42B05, 42C10}
\thanks{\sc Keywords: Iterated Ito stochastic integral,
Iterated stochastic integral with respect to martingale Poisson measures,
Iterated stochastic integral with respect to martingales, 
Generalized Multiple Fourier series, Multiple Fourier--Legendre series,
Mean-square approximation, Expansion.}

\vspace{5mm}

\maketitle {\small
\begin{quote}
\noindent{\sc Abstract.} 
We consider some versions and generalizations 
of the approach
to the expansion of iterated Ito stochastic integrals of 
arbitrary multiplicity $k$ $(k\in\mathbb{N})$ based
on generalized multiple Fourier series.
Expansions 
of iterated stochastic integrals
with respect to martingale Poisson measures 
and with respect to martingales were obtained.
For the iterated stochastic integrals with respect to martingales,
we have proved a theorem which gives a generalization of 
the expansion for iterated Ito stochastic integrals of arbitrary multiplicity 
based
on generalized multiple Fourier series. Also we consider 
a modification of the mentioned expansion of iterated Ito 
stochastic integrals for the case
of complete orthonormal with weight $r(t_1)\ldots r(t_k)\ge 0$ 
systems of functions
in the space $L_2([t, T]^k)$.
Mean-square convergence of the considered expansions is proved.
An example of the expansion of iterated (double) stochastic integrals
with respect to martingales using the system of Bessel functions
is considered.
\medskip
\end{quote}
}

\vspace{12mm}

%\linespread{1.6}

\setlength{\baselineskip}{2.0em}

\tableofcontents

\setlength{\baselineskip}{1.2em}

%\linespread{1.0}

\section{Introduction}

\vspace{5mm}

Let $(\Omega,$ ${\rm F},$ ${\sf P})$ be a complete probability space, let 
$\{{\rm F}_t, t\in[0,T]\}$ be a non-decreasing 
right-continous family of $\sigma$-algebras of ${\rm F},$
and let ${\bf f}_t$ be a standard $m$-dimensional Wiener 
stochastic process which is
${\rm F}_t$-measurable for any $t\in[0, T].$ We assume that the components
${\bf f}_{t}^{(i)}$ $(i=1,\ldots,m)$ of this process are independent. Consider
an Ito stochastic differential equation (SDE) in the integral form

\begin{equation}
\label{1.5.2}
{\bf x}_t={\bf x}_0+\int\limits_0^t {\bf a}({\bf x}_{\tau},\tau)d\tau+
\int\limits_0^t B({\bf x}_{\tau},\tau)d{\bf f}_{\tau},\ \ \
{\bf x}_0={\bf x}(0,\omega).
\end{equation}

\vspace{3mm}
\noindent
Here ${\bf x}_t$ is some $n$-dimensional stochastic process 
satisfying the equation (\ref{1.5.2}). 
The non-random functions ${\bf a}: \mathbb{R}^n\times[0, T]\to
\mathbb{R}^n$,
$B: \mathbb{R}^n\times[0, T]\to\mathbb{R}^{n\times m}$
guarantee the existence and uniqueness up to stochastic equivalence 
of a solution
of the equation (\ref{1.5.2}) \cite{1}. The second integral on the 
right-hand side of (\ref{1.5.2}) is 
interpreted as an Ito stochastic integral.
Let ${\bf x}_0$ be an $n$-dimensional random variable which is 
${\rm F}_0$-measurable and 
${\sf M}\left\{\left|{\bf x}_0\right|^2\right\}<\infty$ 
(${\sf M}$ denotes a mathematical expectation).
We assume that
${\bf x}_0$ and ${\bf f}_t-{\bf f}_0$ are independent when $t>0.$

It is well known \cite{KlPl2}-\cite{KPS}
that Ito SDEs are 
adequate mathematical models of dynamic systems under 
the influence of random disturbances. One of the effective approaches 
to the numerical integration of 
Ito SDEs is an approach based on 
the Taylor--Ito and 
Taylor--Stratonovich expansions
\cite{KlPl2}-\cite{2018aaa}. 
The most important feature of such 
expansions is a presence in them of the so-called iterated
Ito and Stratonovich stochastic integrals which play the key 
role for solving the 
problem of numerical integration of Ito SDEs
and have the 
following form

\begin{equation}
\label{ito}
J[\psi^{(k)}]_{T,t}=\int\limits_t^T\psi_k(t_k) \ldots \int\limits_t^{t_{2}}
\psi_1(t_1) d{\bf w}_{t_1}^{(i_1)}\ldots
d{\bf w}_{t_k}^{(i_k)},
\end{equation}

\vspace{1mm}
\begin{equation}
\label{str}
J^{*}[\psi^{(k)}]_{T,t}=
{\int\limits_t^{*}}^T
\psi_k(t_k) \ldots 
{\int\limits_t^{*}}^{t_2}
\psi_1(t_1) d{\bf w}_{t_1}^{(i_1)}\ldots
d{\bf w}_{t_k}^{(i_k)},
\end{equation}

\vspace{4mm}
\noindent
where every $\psi_l(\tau)$ $(l=1,\ldots,k)$ is a non-random
function 
on $[t,T],$ ${\bf w}_{\tau}^{(i)}={\bf f}_{\tau}^{(i)}$
for $i=1,\ldots,m$ and
${\bf w}_{\tau}^{(0)}=\tau,$\ \
$i_1,\ldots,i_k = 0, 1,\ldots,m,$

\vspace{-1mm}
$$
\int\limits\ \hbox{and}\ \int\limits^{*}
$$ 

\vspace{3mm}
\noindent
denote Ito and 
Stratonovich stochastic integrals,
respectively (in (\ref{str}), 
we use the definition of the Stratonovich stochastic integral from \cite{KlPl2}).

Note that $\psi_l(\tau)\equiv 1$ $(l=1,\ldots,k)$ and
$i_1,\ldots,i_k = 0, 1,\ldots,m$ in  
\cite{KlPl2}-\cite{KlPl1}. At the same time 
$\psi_l(\tau)\equiv (t-\tau)^{q_l}$ ($l=1,\ldots,k,$\ \
$q_1,\ldots,q_k=0, 1, 2,\ldots $) and  $i_1,\ldots,i_k = 1,\ldots,m$ in
\cite{kk5}-\cite{2018aaa}.

The problem of effective jointly numerical modeling 
(with respect to the mean-square convergence criterion) of iterated 
Ito and Stratonovich stochastic integrals 
(\ref{ito}) and (\ref{str}) is 
difficult from 
theoretical and computing point of view \cite{KlPl2}-\cite{KPS},
\cite{2006}-\cite{rr}.

The only exception is connected with a narrow particular case, when 
$i_1=\ldots=i_k\ne 0$ and
$\psi_1(\tau),\ldots,\psi_k(\tau)\equiv \psi(\tau)$.
This case can be investigated
using the Ito formula 
\cite{KlPl2}-\cite{Mi3}.

Note that even for the mentioned coincidence ($i_1=\ldots=i_k\ne 0$)
but for different 
functions $\psi_1(\tau),\ldots,\psi_k(\tau)$ the mentioned 
difficulties persist. As a result, 
relatively simple families of 
iterated Ito and Stratonovich stochastic integrals,
which can be often 
met in the applications, can not be represented effectively in a finite 
form (with respect to 
the mean-square criterion of approximation) using the system of standard 
Gaussian random variables.

Usually, approaches to the expansion of iterated 
stochastic integrals (\ref{ito}) and (\ref{str}) are based on the
expansion of the Wiener process.

For example, in \cite{Mi2} (also see \cite{KlPl2}, \cite{Mi3})
Milstein G.N. 
proposed to expand (\ref{ito}) or 
(\ref{str}) (the case $k=2$ and $i_1\ne i_2;$
$i_1, i_2=1,\ldots,m$)
into the iterated series of products
of standard Gaussian random variables by representing the Brownian
bridge
process as a trigonometric Fourier series with random coefficients 
(the version of the so-called Karhunen--Loeve expansion).
To obtain the Milstein expansion of (\ref{ito}) or
(\ref{str}), the truncated Fourier
expansions of components of the Wiener process ${\bf f}_{\tau}$ must be
iteratively substituted in the single integrals, and the integrals
must be calculated, starting from the innermost integral.
The above procedure
leads to iterated application of 
the operation of limit transition and 
does not lead to a general
expansion of (\ref{ito}) or (\ref{str}) which is valid for an
arbitrary multiplicity $k.$
For this reason, only expansions of single, double, and triple
stochastic integrals were presented 
in \cite{KlPl2} ($k=1, 2, 3$)
and in \cite{Mi2}, \cite{Mi3} ($k=1, 2$) 
for the simplest case $\psi_1(\tau), \psi_2(\tau), \psi_3(\tau)\equiv 1;$ 
$i_1, i_2, i_3=0,1,\ldots,m.$
Moreover, generally speaking, the convergence of approximations 
to the appropriate 
stochastic integrals (\ref{str}) is not proved rigorously
for $k=3$ in \cite{KlPl2}
(Sect.~5.8, pp.~202--204), \cite{KPS} (pp.~82--84),
\cite{KPW} (pp.~438--439),  
\cite{Zapad-9} (pp.~263--264)
(see \cite{2018a}-\cite{12aa-afterxxx} (Sect.~6.2),  \cite{arxiv-1},
\cite{arxiv-3}-\cite{arxiv-7} for details).

Note that in \cite{rrr}, \cite{rr} a method for the expansion
of double \cite{rrr}, \cite{rr} and triple \cite{rrr}
Ito stochastic integrals (\ref{ito}) 
($k=2, 3;$ $\psi_1(\tau),\psi_2(\tau), \psi_3(\tau)$ $\equiv 1;$ $i_1,i_2, 
i_3 =0, 1,\ldots,m$) 
based on the expansion
of the Wiener process using Haar functions \cite{rr} and 
trigonometric functions \cite{rrr}, \cite{rr} has been considered.
The restrictions of this method \cite{rrr}, \cite{rr} are also connected
with the iterated application of the operation of
limit transition at least starting from the second or third 
multiplicity of iterated stochastic integrals.

A more effective and general approach to 
the expansion of iterated Ito stochastic integrals 
(\ref{ito}) of arbitrary multiplicity $k$ $(k\in\mathbb{N})$ based
on generalized multiple Fourier series (converging in the sense of norm
in Hilbert space $L_2([t,T]^k)$)
was proposed and developed by the author of this paper in \cite{2006} (2006)
(also see \cite{2011-2}-\cite{new-2023a},
\cite{301a}-\cite{arxiv-12},
\cite{arxiv-24}-\cite{OK}). 
Hereinafter, this method is referred to as the method of generalized
multiple Fourier series.
As it turned out, the 
method of generalized
multiple Fourier series
can be adapted for the iterated
Stratonovich stochastic integrals (\ref{str}) at least
for the multiplicities 1 to 6 \cite{2011-2}-\cite{12aa-afterxxx}, 
\cite{2010-2}-\cite{2013}, \cite{30a}, \cite{300a},
\cite{400a}, \cite{271a},
\cite{arxiv-5}-\cite{arxiv-8}, \cite{arxiv-23}, \cite{arxiv-9}, \cite{new-art-1-xxy},
\cite{new-art-1xxys}. 
Expansions of these iterated Stratonovich 
stochastic integrals turned out
simpler than the appropriate expansions for
the iterated Ito stochastic integrals (\ref{ito}).

The problem of iterated application of
the operation of limit transition (see above) not appears 
in the method of generalized multiple Fourier series
\cite{2006}-\cite{new-2023a},
\cite{301a}-\cite{arxiv-12},
\cite{arxiv-24}-\cite{OK}.
The idea of this method is as follows: 
the iterated Ito stochastic 
integral (\ref{ito}) of multiplicity $k$ is represented as 
the multiple stochastic 
integral from the certain discontinuous non-random function of $k$ variables
defined on the hypercube $[t, T]^k$, where $[t, T]$ is the interval of 
integration of the iterated Ito stochastic 
integral (\ref{ito}). Then, 
the indicated 
non-random function is expanded in the hypercube into the generalized 
multiple Fourier series converging 
in the mean-square sense
in the space 
$L_2([t,T]^k)$. After a number of nontrivial transformations we come 
(see Theorems 1, 2 below) to the 
mean-square convergening expansion of 
the iterated Ito stochastic 
integral (\ref{ito})
into the multiple 
series of products
of standard  Gaussian random 
variables. The coefficients of this 
series are the coefficients of 
the generalized multiple Fourier series for the mentioned non-random function 
of $k$ variables which can be calculated using the explicit formula 
regardless of the multiplicity $k$ of the
iterated Ito stochastic 
integral (\ref{ito}).
Recall that this method is referred to as the method of generalized
multiple Fourier series.

Thus, we obtain the following new and useful possibilities
of the method of generalized multiple Fourier series.

1. There is the explicit formula (see (\ref{ppppa})) for calculation 
of expansion coefficients 
of the iterated Ito stochastic integral (\ref{ito}) with any
fixed multiplicity $k$. 

2. We have new possibilities for exact calculation of the mean-square 
approximation error
of the iterated Ito stochastic integral (\ref{ito})
of arbitrary
multiplicity $k$
\cite{2017}-\cite{12aa-afterxxx}, \cite{17a}, \cite{arxiv-2}.

3. Since the used
multiple Fourier series is a generalized in the sense
that it is constructed using various complete orthonormal
systems of functions in the space $L_2([t, T])$, then we 
have new possibilities 
for approximation --- we can
use not only the trigonometric functions as in \cite{KlPl2}-\cite{Mi3}
but the Legendre polynomials.

4. As it turned out \cite{2006}-\cite{new-2023a},
\cite{301a}-\cite{arxiv-12},
\cite{arxiv-24}-\cite{OK} it is more convenient to work 
with the Legendre polynomials for constructing the approximations 
of the iterated Ito stochastic integrals (\ref{ito}). 
Approximations based on the Legendre polynomials essentially simpler 
than their analogues based on the trigonometric functions.
Another advantages of the application of Legendre polynomials 
in the framework of the mentioned problem are considered
in \cite{2018a}-\cite{12aa-afterxxx}, \cite{29a}, \cite{301a}.

5. The approach based on the Karhunen--Loeve expansion
of the Brownian bridge process as well as the approach from \cite{rrr},
\cite{rr}
lead to 
iterated application of the operation of
limit transition (the operation of limit 
transition is implemented only once in Theorems 1, 2 (see below))
starting from  
the second multiplicity (in the general case) 
and third multiplicity (for the case
$\psi_1(\tau), \psi_2(\tau), \psi_3(\tau)\equiv 1;$ 
$i_1, i_2, i_3=1,\ldots,m$)
of iterated Ito stochastic integrals.
Multiple series (the operation of limit transition
is implemented only once) are more convenient 
for approximation than the iterated ones
(iterated application of the operation of limit transition), 
since partial sums of multiple series converge for any possible case of  
convergence to infinity of their upper limits of summation 
(let us denote them as $p_1,\ldots, p_k$). 
For example, 
when $p_1=\ldots=p_k=p\to\infty$. For iterated series, 
the condition $p_1=\ldots=p_k=p\to\infty$ obviously 
does not guarantee the convergence of this series.

However, in 
\cite{KlPl2}
(Sect.~5.8, pp.~202--204), \cite{KPS} (pp.~82--84),
\cite{KPW} (pp.~438--439),  
\cite{Zapad-9} (pp.~263--264) the authors use 
(without rigorous proof)
the condition $p_1=p_2=p_3=p\to\infty$
within the frames of the mentioned approach
based on the Karhunen--Loeve expansion of the Brownian bridge
process \cite{Mi2} together with the Wong--Zakai approximation
\cite{W-Z-1}-\cite{Watanabe} 
(see \cite{2018a}-\cite{12aa-afterxxx} (Sect.~6.2),  \cite{arxiv-1},
\cite{arxiv-3}-\cite{arxiv-7} for details).

The method of generalized multiple Fourier series allows
some generalizations and modifications in several directions.

Recently, the method of generalized multiple Fourier series (see Theorems 1, 2 below)
was applied to the expansion and mean-square
approximation of iterated stochastic integrals 
with respect to the infinite-dimensional $Q$-Wiener process
\cite{2018a}-\cite{12aa-afterxxx} (Chapter 7),
\cite{31a}-\cite{200aaa}. These results can be directly
applied to the construction of 
high-order strong numerical methods for non-commutative 
semilinear stochastic partial differential equations
with non-linear multiplicative trace class noise
\cite{2018a}-\cite{12aa-afterxxx} (Chapter 7),
\cite{31a}-\cite{200aaa}.

In this article, we demonstrate that the
method of generalized multiple Fourier series
is essentially
general and allows some transformations for
other types of iterated stochastic integrals.
We will consider versions of the
method of generalized multiple Fourier series
for iterated stochastic integrals
with respect to martingale Poisson measures and 
for iterated stochastic integrals with respect to martingales.
The mentioned results are 
sufficiently natural according to
general properties of martingales.

In Sect.~2, we formulate Theorem 1 on
expansion of the iterated Ito stochastic integrals 
(\ref{ito})
of arbitrary multiplicity $k$
based on generalized multiple Fourier series (method of generalized
multiple Fourier series) \cite{2006}-\cite{new-2023a},
\cite{301a}-\cite{arxiv-12},
\cite{arxiv-24}-\cite{OK}.
Sect.~3 is devoted 
to a generalization of Theorem 1 for the case of an arbitrary 
complete ortho\-nor\-mal system of functions in the space $L_2([t, T])$
and $\psi_1(\tau),$ $\ldots,\psi_k(\tau)\in L_2([t, T]).$
In Sect.~4, we define the stochastic 
integral with respect to the martingale Poisson
measure and consider some properties of this integral.
Sect.~5 is devoted to a version of Theorem 1
for the iterated stochastic integrals
with respect to martingale Poisson measures. In Sect.~6, we 
consider a generalization of Theorem 1 for the case of 
iterated stochastic integrals with respect to martingales.
Sect.~7 is devoted to versions of Theorems 1, 2
for the case of complete orthonormal with weight $r(t_1)\ldots r(t_k)\ge 0$ 
systems of functions in the space $L_2([t, T]^k).$
In Sect.~8, we consider one modification of theorems from Sect.~6 and 7. 
Sect.~9 is devoted to an example of the application of results
from Sect.~8.

\vspace{2mm}

{\it We will say that the function $f(x):$ $[t, T]\to \mathbb{R}^1$ 
satisfies the condition $(\star )$, if it is 
continuous at the interval $[t, T]$ except may be
for the finite number of points 
of the finite discontinuity as well as it is 
right-continuous 
at the interval $[t, T].$}

\vspace{2mm}

Let us suppose that $\{\phi_j(x)\}_{j=0}^{\infty}$
is a complete orthonormal system of functions in the 
space $L_2([t, T]),$ each function  $\phi_j(x)$ 
of which for finite $j$ satisfies the 
condition $(\star)$.

It is clear that
complete orthonormal systems $\{\phi_j(x)\}_{j=0}^{\infty}$ of 
continuous functions in the 
space $L_2([t, T])$ satisfy the condition 
($\star $).

Let us consider some examples of systems satisfying the condition 
($\star $).

\vspace{2mm}

{\bf Example 1.}\ The system of Legendre polynomials

\vspace{-1mm}
\begin{equation}
\label{4009}
\phi_j(x)=\sqrt{\frac{2j+1}{T-t}}P_j\left(\left(
x-\frac{T+t}{2}\right)\frac{2}{T-t}\right),\ \ \ 
j=0, 1, 2,\ldots,\ \ \ x\in [t, T],
\end{equation}

\vspace{3mm}
\noindent
where $P_j(y),$ $y\in [-1,1]$ is the Legendre polynomial

$$
P_j(y)=\frac{1}{2^j j!} \frac{d^j}{dy^j}\left(y^2-1\right)^j.
$$

\vspace{5mm}

{\bf Example 2.} The system of trigonometric functions

$$
\phi_j(x)=\frac{1}{\sqrt{T-t}}
\left\{
\begin{matrix}
1,\ &\ j=0\cr\cr\cr
\sqrt{2}{\rm sin} \left(2\pi r(x-t)/(T-t)\right),\ &\ j=2r-1\cr\cr\cr
\sqrt{2}{\rm cos} \left(2\pi r(x-t)/(T-t)\right),\ &\ j=2r
\end{matrix},\right.
$$

\vspace{3mm}
\noindent
where $x\in[t, T],$\ \ $r=1, 2,\ldots $

\vspace{2mm}

{\bf Example 3.}\ The system of Haar functions

\vspace{1mm}
$$
\phi_0(x)=\frac{1}{\sqrt{T-t}},\ \ \
\phi_{nj}(x)=\frac{1}{\sqrt{T-t}}\varphi_{nj}\biggl(\frac{x-t}{T-t}\biggr),\ \ \
x\in [t, T],
$$

\vspace{4mm}
\noindent
where
$n=0, 1,\ldots,$\ \  $j=1, 2,\ldots, 2^n,$
and the functions $\varphi_{nj}(x)$ have the following form

\vspace{1mm}
$$
\varphi_{nj}(x)=
\begin{cases}
2^{n/2},\ &\ x\in[(j-1)/2^n,\ (j-1)/2^n+
1/2^{n+1})\cr\cr\cr
-2^{n/2},\ &\ x\in[(j-1)/2^n+1/2^{n+1},\
j/2^n)\cr\cr\cr
0,\ &\ \hbox{otherwise}
\end{cases},
$$

\vspace{4mm}
\noindent
where $n=0, 1,\ldots,$\ \  $j=1, 2,\ldots, 2^n$ 
(we choose the values of Haar functions 
in the points of discontinuity in order they will be 
right-continuous).

\vspace{2mm}

{\bf Example 4.}\ The system of 
Rademacher--Walsh functions

$$
\phi_0(x)=\frac{1}{\sqrt{T-t}},\ 
$$

\vspace{2mm}
$$
\phi_{m_1\ldots m_k}(x)=
\frac{1}{\sqrt{T-t}}\varphi_{m_1}\biggl(\frac{x-t}{T-t}\biggr)
\ldots \varphi_{m_k}
\biggl(\frac{x-t}{T-t}\biggr),\ \ \ x\in [t, T],
$$

\vspace{5mm}
\noindent
where $0<m_1<\ldots<m_k,\ \ m_1,\ldots,m_k=1, 2,\ldots,\ \ k=1, 2,\ldots,$

\vspace{1mm}
$$
\varphi_m(x)=(-1)^{[2^m x]},
$$

\vspace{4mm}
\noindent
$x\in [0, 1]$,\ \ $ m=1, 2,\ldots,$\ \ $[y]$ is an integer part of
a real number $y.$

\vspace{5mm}

\section{Method of Expansion of Iterated Ito Stochastic 
Integrals of Arbitrary
Multiplicity Based 
on Generalized Multiple Fourier Series Converging in the Mean}

\vspace{5mm}

Suppose that every $\psi_l(\tau)$ $(l=1,\ldots,k)$ is a 
non-random function from the space $L_2([t, T])$.
Define the following function on the hypercube $[t, T]^k$

\vspace{-1mm}
\begin{equation}
\label{ppp}
K(t_1,\ldots,t_k)=
\begin{cases}
\psi_1(t_1)\ldots \psi_k(t_k)\ &\hbox{for}\ \ t_1<\ldots<t_k\\
~\\
~\\
0\ &\hbox{otherwise}
\end{cases},\ \ \ \ t_1,\ldots,t_k\in[t, T],\ \ \ \ k\ge 2
\end{equation}

\vspace{5mm}
\noindent
and 
$K(t_1)\equiv\psi_1(t_1)$ for $t_1\in[t, T].$

Suppose that $\{\phi_j(x)\}_{j=0}^{\infty}$
is a complete orthonormal system of functions in the space
$L_2([t, T])$. 
The function $K(t_1,\ldots,t_k)$ belongs to the 
space $L_2([t, T]^k).$
At this situation it is well known that the generalized 
multiple Fourier series 
of $K(t_1,\ldots,t_k)\in L_2([t, T]^k)$ is converging 
to $K(t_1,\ldots,t_k)$ in the hypercube $[t, T]^k$ in 
the mean-square sense, i.e.

\vspace{1mm}
$$
\hbox{\vtop{\offinterlineskip\halign{
\hfil#\hfil\cr
{\rm lim}\cr
$\stackrel{}{{}_{p_1,\ldots,p_k\to \infty}}$\cr
}} }\Biggl\Vert
K(t_1,\ldots,t_k)-
\sum_{j_1=0}^{p_1}\ldots \sum_{j_k=0}^{p_k}
C_{j_k\ldots j_1}\prod_{l=1}^{k} \phi_{j_l}(t_l)\Biggr\Vert_{L_2([t,T]^k)}=0,
$$

\vspace{3mm}
\noindent
where
\begin{equation}
\label{ppppa}
C_{j_k\ldots j_1}=\int\limits_{[t,T]^k}
K(t_1,\ldots,t_k)\prod_{l=1}^{k}\phi_{j_l}(t_l)dt_1\ldots dt_k,
\end{equation}

\vspace{1mm}
$$
\left\Vert f\right\Vert_{L_2([t,T]^k)}=\left(\int\limits_{[t,T]^k}
f^2(t_1,\ldots,t_k)dt_1\ldots dt_k\right)^{1/2}.
$$

\vspace{5mm}

Consider the partition $\{\tau_j\}_{j=0}^N$ of $[t,T]$ such that

\begin{equation}
\label{1111}
t=\tau_0<\ldots <\tau_N=T,\ \ \
\Delta_N=
\hbox{\vtop{\offinterlineskip\halign{
\hfil#\hfil\cr
{\rm max}\cr
$\stackrel{}{{}_{0\le j\le N-1}}$\cr
}} }\Delta\tau_j\to 0\ \ \hbox{if}\ \ N\to \infty,\ \ \
\Delta\tau_j=\tau_{j+1}-\tau_j.
\end{equation}

\vspace{4mm}

{\bf Theorem 1}\ \cite{2006} (2006), 
\cite{2011-2}-\cite{new-2023a}, \cite{301a}-\cite{arxiv-12}, 
\cite{arxiv-24}-\cite{OK}.
{\it Suppose that
every $\psi_l(\tau)$ $(l=1,\ldots, k)$ is a continuous non-random func\-tion on 
$[t, T]$ and
$\{\phi_j(x)\}_{j=0}^{\infty}$ is a complete orthonormal
system of functions in the space 
$L_2([t,T]),$ 
each function of which for finite $j$ satisfies the condition 
$(\star)$ {\rm(}see Sect. {\rm 1)}.
Then

\vspace{-1mm}
$$
J[\psi^{(k)}]_{T,t}\  =\ 
\hbox{\vtop{\offinterlineskip\halign{
\hfil#\hfil\cr
{\rm l.i.m.}\cr
$\stackrel{}{{}_{p_1,\ldots,p_k\to \infty}}$\cr
}} }\sum_{j_1=0}^{p_1}\ldots\sum_{j_k=0}^{p_k}
C_{j_k\ldots j_1}\Biggl(
\prod_{l=1}^k\zeta_{j_l}^{(i_l)}\ -
\Biggr.
$$

\vspace{2mm}
\begin{equation}
\label{tyyy}
-\ \Biggl.
\hbox{\vtop{\offinterlineskip\halign{
\hfil#\hfil\cr
{\rm l.i.m.}\cr
$\stackrel{}{{}_{N\to \infty}}$\cr
}} }\sum_{(l_1,\ldots,l_k)\in {\rm G}_k}
\phi_{j_{1}}(\tau_{l_1})
\Delta{\bf w}_{\tau_{l_1}}^{(i_1)}\ldots
\phi_{j_{k}}(\tau_{l_k})
\Delta{\bf w}_{\tau_{l_k}}^{(i_k)}\Biggr),
\end{equation}

\vspace{6mm}
\noindent
where $J[\psi^{(k)}]_{T,t}$ is defined by {\rm (\ref{ito}),}

\vspace{-1mm}
$$
{\rm G}_k={\rm H}_k\backslash{\rm L}_k,\ \ \
{\rm H}_k=\{(l_1,\ldots,l_k):\ l_1,\ldots,l_k=0,\ 1,\ldots,N-1\},
$$

\vspace{-1mm}
$$
{\rm L}_k=\{(l_1,\ldots,l_k):\ l_1,\ldots,l_k=0,\ 1,\ldots,N-1;\
l_g\ne l_r\ (g\ne r);\ g, r=1,\ldots,k\},
$$

\vspace{3mm}
\noindent
${\rm l.i.m.}$ is a limit in the mean-square sense$,$
$i_1,\ldots,i_k=0,1,\ldots,m,$

\vspace{-1mm}
\begin{equation}
\label{rr23}
\zeta_{j}^{(i)}=
\int\limits_t^T \phi_{j}(\tau) d{\bf w}_{\tau}^{(i)}
\end{equation} 

\vspace{2mm}
\noindent
are independent standard Gaussian random variables
for various
$i$ or $j$ {\rm(}if $i\ne 0${\rm),}
$C_{j_k\ldots j_1}$ is the Fourier coefficient {\rm(\ref{ppppa}),}
$\Delta{\bf w}_{\tau_{j}}^{(i)}=
{\bf w}_{\tau_{j+1}}^{(i)}-{\bf w}_{\tau_{j}}^{(i)}$
$(i=0, 1,\ldots,m),$
$\left\{\tau_{j}\right\}_{j=0}^{N}$ is a partition of
the interval $[t, T]$ which satisfies the condition {\rm (\ref{1111})}.
}

\vspace{2mm}

Let us consider transformed particular cases of Theorem 1 for
$k=1,\ldots,5$ \cite{2006}-\cite{new-2023a}, \cite{301a}-\cite{arxiv-12}, 
\cite{arxiv-24}-\cite{OK} (the cases $k=6$ and 7 can be found
in \cite{2011-2}-\cite{2018aaa}, \cite{arxiv-1})

\vspace{1mm}
\begin{equation}
\label{a1}
J[\psi^{(1)}]_{T,t}
=\hbox{\vtop{\offinterlineskip\halign{
\hfil#\hfil\cr
{\rm l.i.m.}\cr
$\stackrel{}{{}_{p_1\to \infty}}$\cr
}} }\sum_{j_1=0}^{p_1}
C_{j_1}\zeta_{j_1}^{(i_1)},
\end{equation}

\vspace{2mm}
\begin{equation}
\label{a2}
J[\psi^{(2)}]_{T,t}
=\hbox{\vtop{\offinterlineskip\halign{
\hfil#\hfil\cr
{\rm l.i.m.}\cr
$\stackrel{}{{}_{p_1,p_2\to \infty}}$\cr
}} }\sum_{j_1=0}^{p_1}\sum_{j_2=0}^{p_2}
C_{j_2j_1}\Biggl(\zeta_{j_1}^{(i_1)}\zeta_{j_2}^{(i_2)}
-{\bf 1}_{\{i_1=i_2\ne 0\}}
{\bf 1}_{\{j_1=j_2\}}\Biggr),
\end{equation}

\vspace{5mm}
$$
J[\psi^{(3)}]_{T,t}=
\hbox{\vtop{\offinterlineskip\halign{
\hfil#\hfil\cr
{\rm l.i.m.}\cr
$\stackrel{}{{}_{p_1,p_2,p_3\to \infty}}$\cr
}} }\sum_{j_1=0}^{p_1}\sum_{j_2=0}^{p_2}\sum_{j_3=0}^{p_3}
C_{j_3j_2j_1}\Biggl(
\zeta_{j_1}^{(i_1)}\zeta_{j_2}^{(i_2)}\zeta_{j_3}^{(i_3)}
-\Biggr.
$$
\begin{equation}
\label{a3}
\Biggl.-{\bf 1}_{\{i_1=i_2\ne 0\}}
{\bf 1}_{\{j_1=j_2\}}
\zeta_{j_3}^{(i_3)}
-{\bf 1}_{\{i_2=i_3\ne 0\}}
{\bf 1}_{\{j_2=j_3\}}
\zeta_{j_1}^{(i_1)}-
{\bf 1}_{\{i_1=i_3\ne 0\}}
{\bf 1}_{\{j_1=j_3\}}
\zeta_{j_2}^{(i_2)}\Biggr),
\end{equation}

\vspace{6mm}

$$
J[\psi^{(4)}]_{T,t}
=
\hbox{\vtop{\offinterlineskip\halign{
\hfil#\hfil\cr
{\rm l.i.m.}\cr
$\stackrel{}{{}_{p_1,\ldots,p_4\to \infty}}$\cr
}} }\sum_{j_1=0}^{p_1}\ldots \sum_{j_4=0}^{p_4}
C_{j_4 \ldots j_1}\Biggl(
\prod\limits_{l=1}^4 \zeta_{j_l}^{(i_l)}
\Biggr.
-
$$
$$
-
{\bf 1}_{\{i_1=i_2\ne 0\}}
{\bf 1}_{\{j_1=j_2\}}
\zeta_{j_3}^{(i_3)}
\zeta_{j_4}^{(i_4)}
-
{\bf 1}_{\{i_1=i_3\ne 0\}}
{\bf 1}_{\{j_1=j_3\}}
\zeta_{j_2}^{(i_2)}
\zeta_{j_4}^{(i_4)}-
$$
$$
-
{\bf 1}_{\{i_1=i_4\ne 0\}}
{\bf 1}_{\{j_1=j_4\}}
\zeta_{j_2}^{(i_2)}
\zeta_{j_3}^{(i_3)}
-
{\bf 1}_{\{i_2=i_3\ne 0\}}
{\bf 1}_{\{j_2=j_3\}}
\zeta_{j_1}^{(i_1)}
\zeta_{j_4}^{(i_4)}-
$$
$$
-
{\bf 1}_{\{i_2=i_4\ne 0\}}
{\bf 1}_{\{j_2=j_4\}}
\zeta_{j_1}^{(i_1)}
\zeta_{j_3}^{(i_3)}
-
{\bf 1}_{\{i_3=i_4\ne 0\}}
{\bf 1}_{\{j_3=j_4\}}
\zeta_{j_1}^{(i_1)}
\zeta_{j_2}^{(i_2)}+
$$
$$
+
{\bf 1}_{\{i_1=i_2\ne 0\}}
{\bf 1}_{\{j_1=j_2\}}
{\bf 1}_{\{i_3=i_4\ne 0\}}
{\bf 1}_{\{j_3=j_4\}}
+
$$
$$
+
{\bf 1}_{\{i_1=i_3\ne 0\}}
{\bf 1}_{\{j_1=j_3\}}
{\bf 1}_{\{i_2=i_4\ne 0\}}
{\bf 1}_{\{j_2=j_4\}}+
$$
\begin{equation}
\label{a4}
+\Biggl.
{\bf 1}_{\{i_1=i_4\ne 0\}}
{\bf 1}_{\{j_1=j_4\}}
{\bf 1}_{\{i_2=i_3\ne 0\}}
{\bf 1}_{\{j_2=j_3\}}\Biggr),
\end{equation}

\vspace{7mm}

$$
J[\psi^{(5)}]_{T,t}
=\hbox{\vtop{\offinterlineskip\halign{
\hfil#\hfil\cr
{\rm l.i.m.}\cr
$\stackrel{}{{}_{p_1,\ldots,p_5\to \infty}}$\cr
}} }\sum_{j_1=0}^{p_1}\ldots \sum_{j_5=0}^{p_5}
C_{j_5 \ldots  j_1}\Biggl(
\prod\limits_{l=1}^5 \zeta_{j_l}^{(i_l)}
-\Biggr.
$$
$$
-
{\bf 1}_{\{i_1=i_2\ne 0\}}
{\bf 1}_{\{j_1=j_2\}}
\zeta_{j_3}^{(i_3)}
\zeta_{j_4}^{(i_4)}
\zeta_{j_5}^{(i_5)}-
{\bf 1}_{\{i_1=i_3\ne 0\}}
{\bf 1}_{\{j_1=j_3\}}
\zeta_{j_2}^{(i_2)}
\zeta_{j_4}^{(i_4)}
\zeta_{j_5}^{(i_5)}-
$$
$$
-
{\bf 1}_{\{i_1=i_4\ne 0\}}
{\bf 1}_{\{j_1=j_4\}}
\zeta_{j_2}^{(i_2)}
\zeta_{j_3}^{(i_3)}
\zeta_{j_5}^{(i_5)}-
{\bf 1}_{\{i_1=i_5\ne 0\}}
{\bf 1}_{\{j_1=j_5\}}
\zeta_{j_2}^{(i_2)}
\zeta_{j_3}^{(i_3)}
\zeta_{j_4}^{(i_4)}-
$$
$$
-
{\bf 1}_{\{i_2=i_3\ne 0\}}
{\bf 1}_{\{j_2=j_3\}}
\zeta_{j_1}^{(i_1)}
\zeta_{j_4}^{(i_4)}
\zeta_{j_5}^{(i_5)}-
{\bf 1}_{\{i_2=i_4\ne 0\}}
{\bf 1}_{\{j_2=j_4\}}
\zeta_{j_1}^{(i_1)}
\zeta_{j_3}^{(i_3)}
\zeta_{j_5}^{(i_5)}-
$$
$$
-
{\bf 1}_{\{i_2=i_5\ne 0\}}
{\bf 1}_{\{j_2=j_5\}}
\zeta_{j_1}^{(i_1)}
\zeta_{j_3}^{(i_3)}
\zeta_{j_4}^{(i_4)}
-{\bf 1}_{\{i_3=i_4\ne 0\}}
{\bf 1}_{\{j_3=j_4\}}
\zeta_{j_1}^{(i_1)}
\zeta_{j_2}^{(i_2)}
\zeta_{j_5}^{(i_5)}-
$$
$$
-
{\bf 1}_{\{i_3=i_5\ne 0\}}
{\bf 1}_{\{j_3=j_5\}}
\zeta_{j_1}^{(i_1)}
\zeta_{j_2}^{(i_2)}
\zeta_{j_4}^{(i_4)}
-{\bf 1}_{\{i_4=i_5\ne 0\}}
{\bf 1}_{\{j_4=j_5\}}
\zeta_{j_1}^{(i_1)}
\zeta_{j_2}^{(i_2)}
\zeta_{j_3}^{(i_3)}+
$$
$$
+
{\bf 1}_{\{i_1=i_2\ne 0\}}
{\bf 1}_{\{j_1=j_2\}}
{\bf 1}_{\{i_3=i_4\ne 0\}}
{\bf 1}_{\{j_3=j_4\}}\zeta_{j_5}^{(i_5)}+
{\bf 1}_{\{i_1=i_2\ne 0\}}
{\bf 1}_{\{j_1=j_2\}}
{\bf 1}_{\{i_3=i_5\ne 0\}}
{\bf 1}_{\{j_3=j_5\}}\zeta_{j_4}^{(i_4)}+
$$
$$
+
{\bf 1}_{\{i_1=i_2\ne 0\}}
{\bf 1}_{\{j_1=j_2\}}
{\bf 1}_{\{i_4=i_5\ne 0\}}
{\bf 1}_{\{j_4=j_5\}}\zeta_{j_3}^{(i_3)}+
{\bf 1}_{\{i_1=i_3\ne 0\}}
{\bf 1}_{\{j_1=j_3\}}
{\bf 1}_{\{i_2=i_4\ne 0\}}
{\bf 1}_{\{j_2=j_4\}}\zeta_{j_5}^{(i_5)}+
$$
$$
+
{\bf 1}_{\{i_1=i_3\ne 0\}}
{\bf 1}_{\{j_1=j_3\}}
{\bf 1}_{\{i_2=i_5\ne 0\}}
{\bf 1}_{\{j_2=j_5\}}\zeta_{j_4}^{(i_4)}+
{\bf 1}_{\{i_1=i_3\ne 0\}}
{\bf 1}_{\{j_1=j_3\}}
{\bf 1}_{\{i_4=i_5\ne 0\}}
{\bf 1}_{\{j_4=j_5\}}\zeta_{j_2}^{(i_2)}+
$$
$$
+
{\bf 1}_{\{i_1=i_4\ne 0\}}
{\bf 1}_{\{j_1=j_4\}}
{\bf 1}_{\{i_2=i_3\ne 0\}}
{\bf 1}_{\{j_2=j_3\}}\zeta_{j_5}^{(i_5)}+
{\bf 1}_{\{i_1=i_4\ne 0\}}
{\bf 1}_{\{j_1=j_4\}}
{\bf 1}_{\{i_2=i_5\ne 0\}}
{\bf 1}_{\{j_2=j_5\}}\zeta_{j_3}^{(i_3)}+
$$
$$
+
{\bf 1}_{\{i_1=i_4\ne 0\}}
{\bf 1}_{\{j_1=j_4\}}
{\bf 1}_{\{i_3=i_5\ne 0\}}
{\bf 1}_{\{j_3=j_5\}}\zeta_{j_2}^{(i_2)}+
{\bf 1}_{\{i_1=i_5\ne 0\}}
{\bf 1}_{\{j_1=j_5\}}
{\bf 1}_{\{i_2=i_3\ne 0\}}
{\bf 1}_{\{j_2=j_3\}}\zeta_{j_4}^{(i_4)}+
$$
$$
+
{\bf 1}_{\{i_1=i_5\ne 0\}}
{\bf 1}_{\{j_1=j_5\}}
{\bf 1}_{\{i_2=i_4\ne 0\}}
{\bf 1}_{\{j_2=j_4\}}\zeta_{j_3}^{(i_3)}+
{\bf 1}_{\{i_1=i_5\ne 0\}}
{\bf 1}_{\{j_1=j_5\}}
{\bf 1}_{\{i_3=i_4\ne 0\}}
{\bf 1}_{\{j_3=j_4\}}\zeta_{j_2}^{(i_2)}+
$$
$$
+
{\bf 1}_{\{i_2=i_3\ne 0\}}
{\bf 1}_{\{j_2=j_3\}}
{\bf 1}_{\{i_4=i_5\ne 0\}}
{\bf 1}_{\{j_4=j_5\}}\zeta_{j_1}^{(i_1)}+
{\bf 1}_{\{i_2=i_4\ne 0\}}
{\bf 1}_{\{j_2=j_4\}}
{\bf 1}_{\{i_3=i_5\ne 0\}}
{\bf 1}_{\{j_3=j_5\}}\zeta_{j_1}^{(i_1)}+
$$
\begin{equation}
\label{a5}
+\Biggl.
{\bf 1}_{\{i_2=i_5\ne 0\}}
{\bf 1}_{\{j_2=j_5\}}
{\bf 1}_{\{i_3=i_4\ne 0\}}
{\bf 1}_{\{j_3=j_4\}}\zeta_{j_1}^{(i_1)}\Biggr),
\end{equation}

\vspace{5mm}
\noindent
where ${\bf 1}_A$ is the indicator of the set $A$.

The convergence in the mean of degree 
$2n$ $(n\in \mathbb{N})$ is proved for approximations from Theorem 1 in 
\cite{2011-2}-\cite{2013}, \cite{arxiv-1}. 
In \cite{2018a}-\cite{2018aaa}, \cite{arxiv-1}-\cite{arxiv-3},
the convergence with probability 1 (further w.~p.~1) is proved for expansions
of iterated Ito stochastic integrals 
of arbitrary  multiplicity $k$ ($k\in\mathbb{N}$)
from Theorem 1 
for the cases of
Legendre polynomials and trigonometric functions.

As follows from Theorem 1, the expansion (\ref{tyyy}) 
is valid for discontinuous
complete
orthonormal systems  
of functions in $L_2([t,T])$ satisfying the condition $(\star )$.
For example, Theorem 1 is valid
for the system of Haar functions 
as well as for the system of Rademacher--Walsh functions
\cite{2006}-\cite{2013}, \cite{arxiv-1}.

\vspace{5mm}

\section{Generalization of Theorem 1 to the Case of an Arbitrary 
Complete Ortho\-nor\-mal System of Functions in the Space $L_2([t, T])$
and $\psi_1(\tau),$ $\ldots,\psi_k(\tau)\in L_2([t, T])$}

\vspace{5mm}

Consider a generalization of formulas (\ref{a1})--(\ref{a5}) 
for the case of arbitrary multiplicity $k$ of 
the iterated Ito stochastic integrals
(\ref{ito}).
In order to do this, let us
consider the unordered
set $\{1, 2, \ldots, k\}$ 
and separate it into two parts:
the first part consists of $r$ unordered 
pairs (sequence order of these pairs is also unimportant) and the 
second one consists of the 
remaining $k-2r$ numbers.
So, we have

\vspace{-3mm}
\begin{equation}
\label{leto5007}
(\{
\underbrace{\{g_1, g_2\}, \ldots, 
\{g_{2r-1}, g_{2r}\}}_{\small{\hbox{part 1}}}
\},
\{\underbrace{q_1, \ldots, q_{k-2r}}_{\small{\hbox{part 2}}}
\}),
\end{equation}

\vspace{5mm}
\noindent
where 
$\{g_1, g_2, \ldots, 
g_{2r-1}, g_{2r}, q_1, \ldots, q_{k-2r}\}=\{1, 2, \ldots, k\},$
braces   
mean an unordered 
set, and pa\-ren\-the\-ses mean an ordered set.

We will say that (\ref{leto5007}) is a partition 
and consider the sum with respect to all possible
partitions

\vspace{-1mm}
\begin{equation}
\label{leto5008}
\sum_{\stackrel{(\{\{g_1, g_2\}, \ldots, 
\{g_{2r-1}, g_{2r}\}\}, \{q_1, \ldots, q_{k-2r}\})}
{{}_{\{g_1, g_2, \ldots, 
g_{2r-1}, g_{2r}, q_1, \ldots, q_{k-2r}\}=\{1, 2, \ldots, k\}}}}
a_{g_1 g_2, \ldots, 
g_{2r-1} g_{2r}, q_1 \ldots q_{k-2r}},
\end{equation}

\vspace{4mm}
\noindent
where $a_{g_1 g_2, \ldots, 
g_{2r-1} g_{2r}, q_1 \ldots q_{k-2r}}\in\mathbb{R}.$

Below there are several examples of sums in the form (\ref{leto5008})

$$
\sum_{\stackrel{(\{g_1, g_2\})}{{}_{\{g_1, g_2\}=\{1, 2\}}}}
a_{g_1 g_2}=a_{12},
$$

\vspace{3mm}
$$
\sum_{\stackrel{(\{\{g_1, g_2\}, \{g_3, g_4\}\})}
{{}_{\{g_1, g_2, g_3, g_4\}=\{1, 2, 3, 4\}}}}
a_{g_1 g_2 g_3 g_4}=a_{12,34} + a_{13,24} + a_{23,14},
$$

\vspace{3mm}
$$
\sum_{\stackrel{(\{g_1, g_2\}, \{q_1, q_{2}\})}
{{}_{\{g_1, g_2, q_1, q_{2}\}=\{1, 2, 3, 4\}}}}
a_{g_1 g_2, q_1 q_{2}}=a_{12,34}+a_{13,24}+a_{14,23}
+a_{23,14}+a_{24,13}+a_{34,12},
$$

\vspace{3mm}
$$
\sum_{\stackrel{(\{g_1, g_2\}, \{q_1, q_{2}, q_3\})}
{{}_{\{g_1, g_2, q_1, q_{2}, q_3\}=\{1, 2, 3, 4, 5\}}}}
a_{g_1 g_2, q_1 q_{2}q_3}
=a_{12,345}+a_{13,245}+a_{14,235}
+a_{15,234}+a_{23,145}+a_{24,135}+
$$
$$
+a_{25,134}+a_{34,125}+a_{35,124}+a_{45,123},
$$

\vspace{3mm}
$$
\sum_{\stackrel{(\{\{g_1, g_2\}, \{g_3, g_{4}\}\}, \{q_1\})}
{{}_{\{g_1, g_2, g_3, g_{4}, q_1\}=\{1, 2, 3, 4, 5\}}}}
a_{g_1 g_2, g_3 g_{4},q_1}
=
a_{12,34,5}+a_{13,24,5}+a_{14,23,5}+
a_{12,35,4}+a_{13,25,4}+a_{15,23,4}+
$$
$$
+a_{12,54,3}+a_{15,24,3}+a_{14,25,3}+a_{15,34,2}+a_{13,54,2}+a_{14,53,2}+
a_{52,34,1}+a_{53,24,1}+a_{54,23,1}.
$$

\vspace{8mm}

Let us consider a generalization of Theorem 1 for the case
of an arbitrary complete orthonormal systems  
of functions in the space $L_2([t,T])$ 
and $\psi_1(\tau),\ldots,\psi_k(\tau)\in L_2([t, T]).$

\vspace{2mm}

{\bf Theorem~2}\ \cite{2018a} (Sect.~1.11), \cite{new-2023a}, \cite{arxiv-1} (Sect.~15).
{\it Suppose that
$\psi_1(\tau),\ldots,\psi_k(\tau)\in L_2([t, T])$ and
$\{\phi_j(x)\}_{j=0}^{\infty}$ is an arbitrary complete orthonormal system  
of functions in the space $L_2([t,T]).$
Then the following expansion

\vspace{1mm}
$$
J[\psi^{(k)}]_{T,t}=
\hbox{\vtop{\offinterlineskip\halign{
\hfil#\hfil\cr
{\rm l.i.m.}\cr
$\stackrel{}{{}_{p_1,\ldots,p_k\to \infty}}$\cr
}} }
\sum\limits_{j_1=0}^{p_1}\ldots
\sum\limits_{j_k=0}^{p_k}
C_{j_k\ldots j_1}\Biggl(
\prod_{l=1}^k\zeta_{j_l}^{(i_l)}+\sum\limits_{r=1}^{[k/2]}
(-1)^r \times
\Biggr.
$$

\vspace{2mm}
\begin{equation}
\label{leto6000}
\times
\sum_{\stackrel{(\{\{g_1, g_2\}, \ldots, 
\{g_{2r-1}, g_{2r}\}\}, \{q_1, \ldots, q_{k-2r}\})}
{{}_{\{g_1, g_2, \ldots, 
g_{2r-1}, g_{2r}, q_1, \ldots, q_{k-2r}\}=\{1, 2, \ldots, k\}}}}
\prod\limits_{s=1}^r
{\bf 1}_{\{i_{g_{{}_{2s-1}}}=~i_{g_{{}_{2s}}}\ne 0\}}
\Biggl.{\bf 1}_{\{j_{g_{{}_{2s-1}}}=~j_{g_{{}_{2s}}}\}}
\prod_{l=1}^{k-2r}\zeta_{j_{q_l}}^{(i_{q_l})}\Biggr)
\end{equation}

\vspace{5mm}
\noindent
con\-verg\-ing in the mean-square sense is valid,
where $[x]$ is an integer part of a real number $x,$
$\prod\limits_{\emptyset}
\stackrel{\sf def}{=}1,$ $\sum\limits_{\emptyset}
\stackrel{\sf def}{=}0;$
another notations are the same as in Theorem~{\rm 1}.}

\vspace{2mm}

In particular from (\ref{leto6000}) for $k=5$ we obtain

\vspace{2mm}

$$
J[\psi^{(5)}]_{T,t}=
\hbox{\vtop{\offinterlineskip\halign{
\hfil#\hfil\cr
{\rm l.i.m.}\cr
$\stackrel{}{{}_{p_1,\ldots,p_5\to \infty}}$\cr
}} }
\sum\limits_{j_1=0}^{p_1}\ldots
\sum\limits_{j_5=0}^{p_5}
C_{j_5 \ldots j_1}
\Biggl(
\prod\limits_{l=1}^5 \zeta_{j_l}^{(i_l)}
-\Biggr.
$$

\vspace{2mm}
$$
-\sum\limits_{\stackrel{(\{g_1, g_2\}, \{q_1, q_{2}, q_3\})}
{{}_{\{g_1, g_2, q_{1}, q_{2}, q_3\}=\{1, 2, 3, 4, 5\}}}}
{\bf 1}_{\{i_{g_{{}_{1}}}=~i_{g_{{}_{2}}}\ne 0\}}
{\bf 1}_{\{j_{g_{{}_{1}}}=~j_{g_{{}_{2}}}\}}
\prod_{l=1}^{3}\zeta_{j_{q_l}}^{(i_{q_l})}+
$$

\vspace{2mm}
$$
\Biggl.+
\sum_{\stackrel{(\{\{g_1, g_2\}, 
\{g_{3}, g_{4}\}\}, \{q_1\})}
{{}_{\{g_1, g_2, g_{3}, g_{4}, q_1\}=\{1, 2, 3, 4, 5\}}}}
{\bf 1}_{\{i_{g_{{}_{1}}}=~i_{g_{{}_{2}}}\ne 0\}}
{\bf 1}_{\{j_{g_{{}_{1}}}=~j_{g_{{}_{2}}}\}}
\Biggl.{\bf 1}_{\{i_{g_{{}_{3}}}=~i_{g_{{}_{4}}}\ne 0\}}
{\bf 1}_{\{j_{g_{{}_{3}}}=~j_{g_{{}_{4}}}\}}
\zeta_{j_{q_1}}^{(i_{q_1})}\Biggr).
$$

\vspace{6mm}

The last equality obviously agrees with
(\ref{a5}).

It should be noted that an analogue of Theorem 2 (the case $i_1,\ldots,i_k=1,\ldots,m$)
was considered 
in \cite{Rybakov1000} using the Hermite polynomials and Wick product.
Note that we use another notations 
\cite{2018a} (Sect.~1.11), \cite{new-2023a}, \cite{arxiv-1} (Sect.~15)
in comparison with \cite{Rybakov1000}.
Moreover, the proof from \cite{Rybakov1000} is different from the proof given in 
\cite{2018a} (Sect.~1.11), \cite{new-2023a}, \cite{arxiv-1} (Sect.~15).
See Sect.~4 in \cite{new-2023a} for details.

Below we demonstrate that an approach
to the expansion of iterated Ito stochastic integrals
considered in Theorems 1, 2 is essentially
general and allows some transformations for
other types of iterated stochastic integrals.

Note that Theorems 1, 2 allow to calculate exactly
the mean-square approximation error of the iterated Ito stochastic
integrals (\ref{ito}) of arbitrary multiplicity $k$
(see \cite{2017-1}--\cite{12aa-afterxxx}, \cite{arxiv-2}).
In these papers 
we consider approxinations
of iterated Ito stochastic
integrals
as the expression on the right-hand side of (\ref{leto6000}) before passing to the limit
with respect to $p_1,\ldots,p_k.$

\vspace{5mm}

\section{Stochastic Integral with Respect to Martingale Poisson Measure}

\vspace{5mm}

Let us consider the Poisson random measure in the space
$[0,T]\times{\bf Y}$
$(\mathbb{R}^n\stackrel{\rm def}{=}{\bf Y})$.
We will denote the values of this measure at the set 
$\Delta\times A$ ($\Delta\subseteq
[0,T],$ $A\subset{\bf Y}$) 
as 
$\nu(\Delta,A).$ Let us assume that

$$
{\sf M}\biggl\{\nu(\Delta,A)\biggr\}=
|\Delta|\Pi(A),
$$ 

\vspace{3mm}
\noindent
where $|\Delta|$ is the Lebesgue measure
of $\Delta,$ 
$\Pi(A)$ is a measure on $\sigma$-algebra $\mathcal{B}$ of Borel
sets of ${\bf Y},$ and $\mathcal{B}_0$ is a subalgebra of $\mathcal{B}$
consisting of sets
$A\subset \mathcal{B}$ which satisfy the condition 
$\Pi(A)<\infty.$

Let us consider the martingale Poisson measure

$$
\tilde\nu(\Delta,A)=\nu(\Delta,A)-|\Delta|\Pi(A).
$$

\vspace{3mm}

Let $(\Omega, {\rm F},{\sf P})$ be a fixed probability
space, let $\{{\rm F}_t,$ $t\in[0,T]\}$
be a non-decreasing family of 
$\sigma$-algebras ${\rm F}_t\subset{\rm F}$.

Assume that:

1. The random variables $\nu([0,t),A)$ are ${\rm F}_t$-measurable
for all $A\subseteq \mathcal{B}_0.$

2. The random variables $\nu([t,t+h),A),$ $A\subseteq \mathcal{B}_0,$
$h>0$ do not depend on $\sigma$-algebra ${\rm F}_t.$

Let us define the class $H_l(\Pi,[0,T])$ 
of random functions  
$\varphi:$ $[0,T]\times{\bf Y}\times\Omega\to\mathbb{R}^1,$
that are ${\rm F}_t$-measurable for all
$t\in[0,T],$ ${\bf y}\in{\bf Y}$ 
and satisfy the
following condition 

\vspace{-1mm}
$$
\int\limits_0^T\int\limits_{{\bf Y}}
{\sf M}\biggl\{|\varphi(t,{\bf y})|^l\biggr\}\Pi(d{\bf y})dt<\infty.
$$

\vspace{2mm}

Let us consider the partition $\{\tau_j\}_{j=0}^N$ of the 
interval $[0,T]$
which 
satisfies
the condition (\ref{1111}).

For $\varphi(t,{\bf y})\in H_2(\Pi,[0,T])$ let us define the stochastic 
integral with respect to the martingale Poisson measure as the following 
mean-square limit 
\cite{1}

\vspace{-1mm}
\begin{equation}
\label{1.10}
\int\limits_0^T\int\limits_{{\bf Y}}
\varphi(t,{\bf y})\tilde\nu(dt,d{\bf y})
\stackrel{\rm def}{=}
\hbox{\vtop{\offinterlineskip\halign{
\hfil#\hfil\cr
{\rm l.i.m.}\cr
$\stackrel{}{{}_{N\to\infty}}$\cr
}} }\int\limits_0^T\int\limits_{{\bf Y}}
\varphi^{(N)}(t,{\bf y})\tilde\nu(dt,d{\bf y}),
\end{equation}

\vspace{2mm}
\noindent
where $\varphi^{(N)}(t,{\bf y})$ is any
sequense of step functions 
from the class $H_2(\Pi,[0,T])$ such that 

\vspace{-1mm}
$$
\hbox{\vtop{\offinterlineskip\halign{
\hfil#\hfil\cr
{\rm lim}\cr
$\stackrel{}{{}_{N\to\infty}}$\cr
}} }
\int\limits_0^T\int\limits_{{\bf Y}}{\sf M}\biggl\{\left|\varphi(t,{\bf y})-
\varphi^{(N)}(t,{\bf y})\right|^2\biggr\}\Pi(d{\bf y})dt= 0.
$$

\vspace{2mm}

It is well known \cite{1} that the stochastic integral 
(\ref{1.10}) exists, it does not depend on selection of the
sequence $\varphi^{(N)}(t,{\bf y})$ and it 
satisfies
w.~p.~1 
the following properties

$$
{\sf M}\left\{\int\limits_0^T
\int\limits_{{\bf Y}}\varphi(t,{\bf y})\tilde\nu(dt,d{\bf y})
\biggl|\biggr.{\rm F}_0\right\}=0,
$$

\vspace{1mm}
$$
\int\limits_0^T\int\limits_{{\bf Y}}
(\alpha\varphi_1(t,{\bf y})+\beta\varphi_2(t,{\bf y}))
\tilde\nu(dt,d{\bf y})=
\alpha\int\limits_0^T
\int\limits_{{\bf Y}}\varphi_1(t,{\bf y})\tilde\nu(dt,d{\bf y})
+\beta\int\limits_0^T\int\limits_{{\bf Y}}
\varphi_2(t,{\bf y})\tilde\nu(dt,d{\bf y}),
$$

\vspace{3mm}
$$
{\sf M}\left\{\left|\int\limits_0^T
\int\limits_{{\bf Y}}\varphi(t,{\bf y})\tilde\nu(dt,d{\bf y})
\right|^2 \biggl|\biggr.{\rm F}_0\right\}=
\int\limits_0^T\int\limits_{{\bf Y}}{\sf M}\biggl\{\left|\varphi(t,{\bf y})
\right|^2 \biggl|\biggr.{\rm F}_0
\biggr\}\Pi(d{\bf y})dt,
$$

\vspace{7mm}
\noindent
where $\alpha,$ $\beta$ are some real constants and $\varphi_1(t,{\bf y}),$
$\varphi_2(t,{\bf y}),$
$\varphi(t,{\bf y})$ from the class 
$H_2(\Pi,[0,T]).$

The stochastic integral

\vspace{-1mm}
$$
\int\limits_0^T\int\limits_{{\bf Y}}
\varphi(t,{\bf y})\nu(dt,d{\bf y})
$$

\vspace{3mm}
\noindent
with respect to the Poisson random measure will be defined as follows \cite{1}

\vspace{-1mm}
$$
\int\limits_0^T\int\limits_{{\bf Y}}
\varphi(t,{\bf y})\nu(dt,d{\bf y})=
\int\limits_0^T\int\limits_{{\bf Y}}
\varphi(t,{\bf y})\tilde\nu(dt,d{\bf y})+
\int\limits_0^T\int\limits_{{\bf Y}}
\varphi(t,{\bf y})\Pi(d{\bf y})dt,
$$

\vspace{3mm}
\noindent
where we suppose that the right-hand side of the last equality
exists.

According to the Ito formula for Ito processes with 
jump 
component, we obtain w.~p.~1
\cite{1}

\vspace{-1mm}
\begin{equation}
\label{16.008}
\left({z}_t\right)^n=\int\limits_0^t\int\limits_{{\bf Y}}
\biggl(({z}_{\tau-}+\gamma(\tau,{\bf y}))^n-\left({z}_{\tau-}\right)^n\biggr)
\nu(d\tau,d{\bf y}),
\end{equation}

\vspace{2mm}
\noindent
where $n\in\mathbb{N},$
$$
{z}_t=\int\limits_0^t\int\limits_{{\bf Y}}
\gamma(\tau,{\bf y})
\nu(d\tau,d{\bf y}).
$$

\vspace{3mm}

We suppose that the function $\gamma(\tau,{\bf y})$
satisfies
the conditions of existence of the right-hand side 
of (\ref{16.008}) \cite{1}.

Let us consider \cite{1} the useful estimate for moments
of the stochastic 
integral with respect to the Poisson random measure

\vspace{-1mm}
\begin{equation}
\label{16.010}
a_n(T)\le
\max\limits_{j\in\{n,\ 1\}}
\left\{\left(\int\limits_0^T\int\limits_{{\bf Y}}\left(
\left(\left(b_n(\tau,{\bf y})\right)^{1/n}+1\right)^n-1\right)
\Pi(d{\bf y})d\tau\right)^j\right\},
\end{equation}

\vspace{2mm}
\noindent
where 
$$
a_n(t)=\sup\limits_{0\le\tau\le t}
{\sf M}\biggl\{|{z}_{\tau}|^n\biggr\},\ \ \
b_n(\tau,{\bf y})={\sf M}\biggl\{\left|
\gamma(\tau,{\bf y})\right|^n\biggr\}.
$$

\vspace{7mm}

We suppose that the right-hand side 
of 
(\ref{16.010}) exists. Since

\vspace{1mm}
$$
\tilde\nu(dt,d{\bf y})=\nu(dt,d{\bf y})-\Pi(d{\bf y})dt,
$$

\vspace{4mm}
\noindent
then
according to the Minkowski inequality, we obtain

\vspace{1mm}
\begin{equation}
\label{16.011}
\left({\sf M}\biggl\{\left|\tilde{z}_{t}\right|^{2n}\biggr\}
\right)^{1/2n}\le
\left({\sf M}\biggl\{\left|{z}_{t}\right|^{2n}\biggr\}\right)^{1/2n}+
\left({\sf M}\biggl\{\left|\hat{z}_{t}\right|^{2n}\biggr\}
\right)^{1/2n},
\end{equation}

\vspace{3mm}
\noindent
where
$$
\hat{z}_{t}\stackrel{\rm def}{=}\int\limits_0^t\int\limits_{{\bf Y}}
\gamma(\tau,{\bf y})\Pi(d{\bf y})d\tau
$$
and
$$
\tilde{z}_t=\int\limits_0^t\int\limits_{{\bf Y}}
\gamma(\tau,{\bf y})
\tilde\nu(d\tau,d{\bf y}).
$$

\vspace{2mm}

The value ${\sf M}\biggl\{|\hat{z}_{\tau}|^{2n}\biggr\}$
can be estimated using the 
well known inequality \cite{1}

\begin{equation}
\label{dur}
{\sf M}\left\{|\hat{z}_t|^{2n}\right\}\le t^{2n-1}
\int\limits_{0}^t 
{\sf M}\left\{\left|\int\limits_{{\bf Y}}
\varphi(\tau,{\bf y})\Pi(d{\bf y})\right|^{2n}\right\}d\tau,
\end{equation}

\vspace{2mm}
\noindent
where we suppose that 

$$
\int\limits_0^t{\sf M}\left\{\left|\int\limits_{{\bf Y}}
\gamma(\tau,{\bf y})\Pi(d{\bf y})\right|^{2n}\right\}d\tau<\infty.
$$

\vspace{5mm}

\section{Expansion of Iterated Stochastic Integrals 
with Respect to Martingale Poisson Measures
Based on Generalized Multiple Fourier Series}

\vspace{5mm} 

Let us consider the following iterated stochastic integrals 

\begin{equation}
\label{2000.2.1}
P[\chi^{(k)}]_{T,t}
=\int\limits_t^T\int\limits_{\bf X}\chi_k(t_k,{\bf y}_k)\ldots
\int\limits_t^{t_2}\int\limits_{\bf X}\chi_1(t_1,{\bf y}_1)
\tilde\nu^{(i_1)}(dt_1,d{\bf y}_1)\ldots
\tilde\nu^{(i_k)}(dt_k,d{\bf y}_k),
\end{equation}

\vspace{3mm}
\noindent
where $i_1,\ldots,i_k=0, 1,\ldots,m,$ 
$\mathbb{R}^n\stackrel{\rm def}{=}{\bf X},$

$$
\chi_l(\tau,{\bf y})=\psi_l(\tau)\varphi_l({\bf y})\ \ \ (l=1,\ldots,k),
$$

\vspace{4mm}
\noindent
every function
$\psi_l(\tau): [t,T]\to \mathbb{R}^1$ $(l=1,\ldots,k)$ and every function 
$\varphi_l({\bf y}): {\bf X}\to \mathbb{R}^1$ $(l=1,\ldots,k)$ is such that

\vspace{-2mm}
$$
\chi_l(s,{\bf y})\in
H_2(\Pi,[t,T])\ \ \ (l=1,\ldots,k),
$$

\vspace{4mm}
\noindent 
where definition of the class $H_2(\Pi,[t,T])$ is given above,

$$
\nu^{(i)}(dt,d{\bf y})\ \ \ (i=1,\ldots,m)
$$ 

\vspace{3mm}
\noindent
are independent Poisson random measures for various $i$
which are defined on
$[0,T]\times {\bf X}$,

$$
\tilde\nu^{(i)}(dt,d{\bf y})=
\nu^{(i)}(dt,d{\bf y})-\Pi(d{\bf y})dt\ \ \ (i=1,\ldots,m)
$$ 

\vspace{3mm}
\noindent
are independent martingale Poisson measures for various $i$,

\vspace{-1mm}
$$
\tilde\nu^{(0)}(dt,d{\bf y})\stackrel{\rm def}{=}\Pi(d{\bf y})dt.
$$

\vspace{3mm}

Let us formulate an analoque of Theorem 1 for the iterated
stochastic integrals (\ref{2000.2.1}).

\vspace{2mm}

{\bf Theorem 3}\ \cite{2017}-\cite{12aa-afterxxx}.\ 
{\it Suppose that the following 
conditions are fulfilled{\rm :}

{\rm 1}.\ Every $\psi_l(\tau)\ (l=1,\ldots,k)$ is a 
continuous non-random function at
the interval $[t, T]$.

{\rm 2}.\ $\{\phi_j(x)\}_{j=0}^{\infty}$ is a complete orthonormal
system of functions in the space 
$L_2([t,T]),$ 
each function of which for finite $j$ satisfies the condition 
$(\star)$ {\rm(}see Sect. {\rm 1)}.

{\rm 3}.\ For $l=1,\ldots,k$ and $q=2^{k+1}$
the following condition is fulfilled

\vspace{-1mm}
$$
\int\limits_{\bf X}\left|\varphi_l({\bf y})\right|^q
\Pi(d{\bf y})<\infty.
$$

\vspace{2mm}

Then, for the iterated stochastic integral with respect to 
martingale Poisson measures $P[\chi^{(k)}]_{T,t}$ defined by
{\rm (\ref{2000.2.1})}
the following expansion

\vspace{1mm}
$$
P[\chi^{(k)}]_{T,t}=
\hbox{\vtop{\offinterlineskip\halign{
\hfil#\hfil\cr
{\rm l.i.m.}\cr
$\stackrel{}{{}_{p_1,\ldots,p_k\to \infty}}$\cr
}} }\sum_{j_1=0}^{p_1}\ldots\sum_{j_k=0}^{p_k}
C_{j_k\ldots j_1}
\Biggl(
\prod_{g=1}^k\pi_{j_g}^{(g,i_g)}
\Biggr.-
$$

\vspace{2mm}
\begin{equation}
\label{tyyys}
-\Biggl.
\hbox{\vtop{\offinterlineskip\halign{
\hfil#\hfil\cr
{\rm l.i.m.}\cr
$\stackrel{}{{}_{N\to \infty}}$\cr
}} }\sum_{(l_1,\ldots,l_k)\in {\rm G}_k}
\prod_{g=1}^k 
\phi_{j_{g}}(\tau_{l_g})
\int\limits_{\bf X}\varphi_g({\bf y})
\tilde \nu^{(i_g)}([\tau_{l_g},\tau_{l_g+1}),d{\bf y})
\Biggr)
\end{equation}

\vspace{5mm}
\noindent
converging in the mean-square sense  
is valid, where $\left\{\tau_{j}\right\}_{j=0}^{N}$ is a partition of
the interval $[t, T]$ which satisfies the condition {\rm (\ref{1111}),}

\vspace{1mm}
$$
{\rm G}_k={\rm H}_k\backslash{\rm L}_k,\ \ \
{\rm H}_k=\{(l_1,\ldots,l_k):\ l_1,\ldots,l_k=0,\ 1,\ldots,N-1\},
$$

\vspace{-1mm}
$$
{\rm L}_k=\{(l_1,\ldots,l_k):\ l_1,\ldots,l_k=0,\ 1,\ldots,N-1;\
l_g\ne l_r\ (g\ne r);\ g, r=1,\ldots,k\},
$$

\vspace{5mm}
\noindent
${\rm l.i.m.}$ is a limit in the mean-square sense$,$
$i_1,\ldots,i_k=0,1,\ldots,m,$ 
random variables

\vspace{-1mm}
$$
\pi_{j}^{(g,i_g)}=
\int\limits_t^T \phi_j(\tau)\int\limits_{\bf X}\varphi_g({\bf y})
\tilde\nu^{(i_g)}(d\tau,d{\bf y})
$$ 

\vspace{2mm}
\noindent
are independent for various
$i_g\ne 0$ and uncorrelated for various
$j,$

\vspace{1mm}
$$
C_{j_k\ldots j_1}=\int\limits_{[t,T]^k}
K(t_1,\ldots,t_k)\prod_{l=1}^{k}\phi_{j_l}(t_l)dt_1\ldots dt_k
$$

\vspace{2mm}
\noindent
is the Fourier coefficient,

$$
K(t_1,\ldots,t_k)=
\begin{cases}
\psi_1(t_1)\ldots \psi_k(t_k)\ &\hbox{for}\ \ t_1<\ldots<t_k\\
~\\
~\\
0\ &\hbox{otherwise}
\end{cases},\ \ \ \ t_1,\ldots,t_k\in[t, T],\ \ \ \ k\ge 2
$$

\vspace{5mm}
\noindent
and 
$K(t_1)\equiv\psi_1(t_1)$ for $t_1\in[t, T].$  
}

\vspace{5mm}

{\bf Proof.} The scheme of the proof of Theorem 3 is the same 
as the scheme
of the proof of Theorem 1 
(see \cite{2006}-\cite{2013}, \cite{arxiv-1} for details).
Some differences will take place in 
the proof of the following lemmas (Lemmas 1, 2) and in the final part 
of the proof of Theorem 3.

\vspace{2mm}

{\bf Lemma 1}\ \cite{2011-2}-\cite{2013}.\ {\it Suppose that
every $\psi_l(\tau)$ $(l=1,\ldots,k)$ is a continuous 
function at the interval
$[t, T]$ and every function $\varphi_l({\bf y})$ $(l=1,\ldots,k)$ is such that

\vspace{-1mm}
$$
\int\limits_{\bf X}\left|\varphi_l({\bf y})\right|^2
\Pi(d{\bf y})<\infty.
$$ 

\vspace{2mm}

Then, the following equality 

\vspace{-1mm}
\begin{equation}
\label{2000.2.11}
P[\bar \chi^{(k)}]_{T,t}=
\hbox{\vtop{\offinterlineskip\halign{
\hfil#\hfil\cr
{\rm l.i.m.}\cr
$\stackrel{}{{}_{N\to \infty}}$\cr
}} }
\sum_{j_k=0}^{N-1}
\ldots \sum_{j_1=0}^{j_{2}-1}
\prod_{l=1}^k\int\limits_{\bf X}\chi_l(\tau_{j_l},{\bf y})
\bar \nu^{(i_l)}([\tau_{j_l},\tau_{j_l+1}),d{\bf y})
\end{equation}

\vspace{3mm}
\noindent
is valid w.~p.~{\rm 1,}
where $\left\{\tau_{j}\right\}_{j=0}^{N}$ is a partition of
the interval $[t,T]$ which satisfies the condition {\rm (\ref{1111})},

\vspace{1mm}
$$
\bar \nu^{(i)}([\tau,s),d{\bf y})=
\begin{cases}
\tilde \nu^{(i)}([\tau,s),d{\bf y})\cr\cr
\nu^{(i)}([\tau,s),d{\bf y})
\end{cases} \ \ \ (i=0, 1, \ldots,m),
$$

\vspace{5mm}
\noindent
the integral $P[\bar \chi^{(k)}]_{T,t}$
differs from the integral $P[\chi^{(k)}]_{T,t}$
{\rm (}see {\rm (\ref{2000.2.1}))} by the fact that in 
$P[\bar \chi^{(k)}]_{T,t}$
we use $\bar\nu^{(i_l)}(dt_l,d{\bf y}_l)$
instead of 
$\tilde\nu^{(i_l)}(dt_l,d{\bf y}_l)$
$(l=1,\ldots,k).$
}

\vspace{2mm}

{\bf Proof.} Using the moment properties of stochastic 
integrals with respect 
to 
Poisson random measures (see above) and conditions of Lemma 1, 
it is easy to notice that the integral sum 
of the integral 
$P[\bar\chi^{(k)}]_{T,t}$ 
under
the conditions of Lemma 1 can be represented as a sum of the 
expression from the right-hand side of (\ref{2000.2.11}) 
before passing to the limit
$\hbox{\vtop{\offinterlineskip\halign{
\hfil#\hfil\cr
{\rm l.i.m.}\cr
$\stackrel{}{{}_{N\to \infty}}$\cr
}} }$
and the 
value which 
converges
to zero 
in the mean-square sense if 
$N\to \infty.$

Note that in the case when the functions 
$\psi_l(\tau)$ $(l=1,\ldots,k)$ satisfy the condition
$(\star)$ {\rm (}see Sect. {\rm 1)}
we can suppose that among the points
$\tau_j,$ $j=0,1,\ldots,N$ there are all points of 
jumps of the functions $\psi_l(\tau)$ 
$(l=1,\ldots,k)$. Further,
we can apply the argumentation as in Sect.~4 from
\cite{arxiv-1} (also see \cite{2006}-\cite{12aa-afterxxx}).

Let us consider the following multiple and iterated
stochastic integrals 

\vspace{2mm}
$$
\hbox{\vtop{\offinterlineskip\halign{
\hfil#\hfil\cr
{\rm l.i.m.}\cr
$\stackrel{}{{}_{N\to \infty}}$\cr
}} }
\sum_{j_1,\ldots,j_k=0}^{N-1}
\Phi(\tau_{j_1},\ldots,\tau_{j_k})
\prod_{l=1}^k 
\int\limits_{\bf X}\varphi_l({\bf y})
\tilde \nu^{(i_l)}([\tau_{j_l},\tau_{j_l+1}),d{\bf y})
\stackrel{\rm def}{=}P[\Phi]_{T,t}^{(k)},
$$

\vspace{3mm}
$$
\int\limits_t^T
\ldots\int\limits_t^{t_2}
\Phi(t_{1},\ldots,t_{k})\int\limits_{\bf X}\varphi_1({\bf y})
\tilde\nu^{(i_1)}(dt_1,d{\bf y})\ldots
\int\limits_{\bf X}\varphi_k({\bf y})\tilde\nu^{(i_k)}(dt_k,d{\bf y})
\stackrel{\rm def}{=}
\hat P[\Phi]_{T,t}^{(k)},
$$

\vspace{6mm}
\noindent
where the sense of notations of the formula 
(\ref{2000.2.11}) is saved and
$\Phi(t_1,\ldots,t_k):$ $[t, T]^k\to\mathbb{R}^1$ is a bounded non-random
function.

Note that if the functions $\varphi_l({\bf y})$ $(l=1,\ldots,k)$
satisfy the conditions of Lemma 1 and the function 
$\Phi(t_1,\ldots,t_k)$
is continuous in the domain of integration, then for the integral 
$\hat P[\Phi]_{T,t}^{(k)}$ the equality similar to
{\rm (\ref{2000.2.11})}
is valid w.~p.~1.

\vspace{2mm}

{\bf Lemma 2}\ \cite{2011-2}-\cite{2013}.\ {\it Assume that
the following conditions 
are fulfilled{\rm :}

$$
g_l(\tau,{\bf y})=h_l(\tau)\varphi_l({\bf y})\ \ \ (l=1,\ldots,k),
$$

\vspace{4mm}
\noindent
where 
the functions $h_l(\tau):$ $[t, T]\to\mathbb{R}^1$ 
$(l=1,\ldots,k)$ satisfy the condition
$(\star)$ {\rm (}see Sect. {\rm 1)} and 
the functions 
$\varphi_l({\bf y}):$ ${\bf X}\to\mathbb{R}^1$ 
$(l=1,\ldots,k)$ satisfy the condition

$$
\int\limits_{\bf X}\left|\varphi_l({\bf y})
\right|^p
\Pi(d{\bf y})<\infty\ \ \ \hbox{for}\ \ \ p=2^{k+1}.
$$

\vspace{2mm}

Then
$$
\prod_{l=1}^k \int\limits_t^T\int\limits_{\bf X} g_l(s,{\bf y}) 
\bar\nu^{(i_l)}(ds,d{\bf y})=
P[\Phi]_{T,t}^{(k)}\ \ \ \hbox{w.~p.~{\rm 1}},
$$

\vspace{4mm}
\noindent
where $i_l=0, 1, \ldots,m$ $(l=1,\ldots,k)$ and
$$
\Phi(t_1,\ldots,t_k)
=\prod\limits_{l=1}^k h_l(t_l).
$$
}

\vspace{2mm}

{\bf Proof.} Let us introduce the following notations

\vspace{2mm}
$$
J[\bar g_l]_N\stackrel{\rm def}{=}\sum\limits_{j=0}^{N-1}\int\limits_{\bf X}
g_l(\tau_j,{\bf y})\bar\nu^{(i_l)}([\tau_{j},\tau_{j+1}),d{\bf y}),
$$

$$
J[\bar g_l]_{T,t}
\stackrel{\rm def}{=}\int\limits_t^T
\int\limits_{\bf X}
g_l(s,{\bf y})
\bar\nu^{(i_l)}(ds,d{\bf y}),
$$

\vspace{5mm}
\noindent
where $\{\tau_j\}_{j=0}^N$ is a partition of the interval 
$[t,T]$ satisfying the condition (\ref{1111}).

It is easy to see that 

$$
\prod_{l=1}^k J[\bar g_l]_N-\prod_{l=1}^k J[\bar g_l]_{T,t}=
$$

\vspace{1mm}
$$
=
\sum_{l=1}^k \left(\prod_{q=1}^{l-1} J[\bar g_q]_{T,t}\right)
\left(J[\bar g_l]_N-
J[\bar g_l]_{T,t}\right)\left(\prod_{q=l+1}^k J[\bar g_q]_N\right).
$$

\vspace{4mm}

Using the Minkowski inequality and the inequality of Cauchy--Bunyakovsky
together with estimates 
of moments of stochastic integrals with respect to Poisson random
measures (see Sect.~4) and conditions 
of Lemma 2, we obtain

\vspace{-1mm}
\begin{equation}
\label{2000.4.3}
\left({\sf M}\left\{\left|\prod_{l=1}^k J[\bar g_l]_N-
\prod_{l=1}^k J[\bar g_l]_{T,t}\right|^2\right\}
\right)^{1/2}\le C_k
\sum_{l=1}^k
\left({\sf M}
\left\{\biggl|J[\bar g_l]_N-J[\bar g_l]_{T,t}
\biggr|^4\right\}\right)^{1/4},
\end{equation}

\vspace{3mm}
\noindent
where $C_k<\infty.$

We have

\vspace{-2mm}
$$
J[\bar g_l]_N-J[\bar g_l]_{T,t}
=\sum\limits_{q=0}^{N-1}J[\Delta\bar g_{l}]_{\tau_{q+1},\tau_q},
$$

\vspace{4mm}
\noindent
where
$$
J[\Delta\bar g_{l}]_{\tau_{q+1},\tau_q}
=\int\limits_{\tau_q}^{\tau_{q+1}}\int\limits_{\bf X}\left(
g_l(\tau_q,{\bf y})-g_l(s,{\bf y})\right)
\bar\nu^{(i_l)}(ds,d{\bf y}).
$$

\vspace{3mm}

Let as introduce the notation

$$
h_l^{(N)}(s)=h_l(\tau_q),\ \ \ s\in [\tau_q, \tau_{q+1}), \ \ \ q=0, 1, \ldots,
N-1.
$$

\vspace{2mm}

Then
$$
J[\bar g_l]_N-J[\bar g_l]_{T,t}
=\sum\limits_{q=0}^{N-1}J[\Delta\bar g_{l}]_{\tau_{q+1},\tau_q}=
$$
$$
=
\int\limits_{t}^{T}\left(h_l^{(N)}(s)-
h_l(s)\right)\int\limits_{\bf X}
\phi_l({\bf y})
\bar\nu^{(i_l)}(ds,d{\bf y}).
$$

\vspace{4mm}

Applying the estimate (\ref{16.010}) for $n=4$ 
and the estimates (\ref{16.011}), (\ref{dur}) for $n=2$ to
the value

$$
{\sf M}\left\{\left|
\int\limits_{t}^{T}\left(h_l^{(N)}(s)-
h_l(s)\right)\int\limits_{\bf X}
\phi_l({\bf y})
\bar\nu^{(i_l)}(ds,d{\bf y})\right|^4\right\},
$$

\vspace{3mm}
\noindent
taking into account (\ref{2000.4.3}) together with 
the conditions of Lemma 2 and
the following estimate

\vspace{-3mm}
\begin{equation}
\label{dur2}
\left|h_l(\tau_q)-
h_l(s)\right|<\varepsilon,\ \ \ s\in[\tau_q,\tau_{q+1}],\ \ \ q=0,1,\ldots,N-1,
\end{equation}

\vspace{2mm}
\noindent
where $\varepsilon$ is an arbitrary small positive real number,
we obtain that the right-hand side of (\ref{2000.4.3}) 
converges
to zero when $N\to\infty.$
Considering this fact, we come to 
the statement of Lemma 2. 

It should be noted that (\ref{dur2}) is valid
if the functions $h_l(s)$ are continuous at the interval
$[t, T]$, i.e. these functions are uniformly continuous at this interval.
So, $\left|h_l(\tau_q)-h_l(s)\right|<\varepsilon$
if $s\in [\tau_q, \tau_{q+1}],$ where
$|\tau_{q+1}-\tau_q|<\delta(\varepsilon),$ $q=0, 1,\ldots,N-1$
($\delta(\varepsilon)>0$ exists
for any $\varepsilon>0$ and it does not
depend on points of the interval $[t, T]$).

In the case when the functions 
$h_l(s)$ 
$(l=1,\ldots,k)$ satisfy the condition
$(\star)$ {\rm (}see Sect.~{\rm 1)}
we can suppose that among the points
$\tau_q,$ $q=0,1,\ldots,N$ there are all points of 
jumps of the functions $h_l(s)$ 
$(l=1,\ldots,k)$. Further,
we can apply the argumentation as in Sect.~4 from
\cite{arxiv-1} (also see \cite{2006}-\cite{12aa-afterxxx}).

Obviously, if $i_l=0$ for some $l=1,\ldots,k,$ then
we also come to the 
statement of Lemma 2. Lemma 2 is proved.

Proving Theorem 3 according to the scheme used for the proof
of Theorem 1 in \cite{arxiv-1} 
or Theorem 1.1 in
\cite{2018a}-\cite{12aa-afterxxx} (also see \cite{2006} (Theorem 5.1, P. 236-237), 
\cite{2017} (Theorem 1, P. A.22-A.23),
\cite{2017-1} (Theorem 5.1, P. A.250), 
\cite{2018} (Theorem 5.1, P. A.252-A.253)) and
using Lemmas 1, 2 together with estimates for moments
of stochastic integrals with respect to Poisson random measures
(see Sect.~4), we obtain

$$
{\sf M}\left\{\biggl(R_{T,t}^{p_1,\ldots,p_k}\biggr)^2\right\}
\le 
$$

\vspace{2mm}
$$
\le
C_k  \prod\limits_{l=1}^k
\int\limits_{\bf X}\varphi_l^2({\bf y})\Pi(d{\bf y})
\sum_{(t_1,\ldots,t_k)}
\int\limits_{t}^{T}
\ldots
\int\limits_{t}^{t_2}
\left(K(t_1,\ldots,t_k)-
\sum_{j_1=0}^{p_1}\ldots
\sum_{j_k=0}^{p_k}
C_{j_k\ldots j_1}
\prod_{l=1}^k\phi_{j_l}(t_l)\right)^2\times
$$

\vspace{1mm}
\begin{equation}
\label{qu222}
\times
dt_1
\ldots
dt_k=
\end{equation}

\vspace{2mm}
$$
=
C_k  \prod\limits_{l=1}^k
\int\limits_{\bf X}\varphi_l^2({\bf y})\Pi(d{\bf y})
\int\limits_{[t,T]^k}
\Biggl(K(t_1,\ldots,t_k)-
\sum_{j_1=0}^{p_1}\ldots
\sum_{j_k=0}^{p_k}
C_{j_k\ldots j_1}
\prod_{l=1}^k\phi_{j_l}(t_l)\Biggr)^2\times
$$

$$
\times
dt_1
\ldots
dt_k\le
$$

\vspace{2mm}
$$
\le \bar{C_k}
\int\limits_{[t,T]^k}
\Biggl(K(t_1,\ldots,t_k)-
\sum_{j_1=0}^{p_1}\ldots
\sum_{j_k=0}^{p_k}
C_{j_k\ldots j_1}
\prod_{l=1}^k\phi_{j_l}(t_l)\Biggr)^2
dt_1
\ldots
dt_k\to 0
$$

\vspace{4mm}
\noindent
if $p_1,\ldots,p_k\to\infty,$
where constant $\bar{C_k}$ depends only on $k$ (multiplicity
of the 
iterated stochastic integral with respect to martingale Poisson measures).
At that permutations $(t_1,\ldots,t_k)$ when summing

\vspace{-1mm}
$$
\sum_{(t_1,\ldots,t_k)}
$$

\vspace{2mm}
\noindent
in (\ref{qu222})
are performed only in the values $dt_1\ldots dt_k$ and
indexes near 
upper limits of integration are changed 
correspondently.
Moreover, $R_{T,t}^{p_1,\ldots,p_k}$ has the following form

\vspace{1mm}
$$
R_{T,t}^{p_1,\ldots,p_k}
=\sum_{(t_1,\ldots,t_k)}
\int\limits_{t}^{T}
\ldots
\int\limits_{t}^{t_2}
\left(K(t_1,\ldots,t_k)-
\sum_{j_1=0}^{p_1}\ldots
\sum_{j_k=0}^{p_k}
C_{j_k\ldots j_1}
\prod_{l=1}^k\phi_{j_l}(t_l)\right)\times
$$

\vspace{2mm}
\begin{equation}
\label{jter}
\times
\int\limits_{\bf X}\varphi_1({\bf y})
\tilde \nu^{(i_1)}(dt_1,d{\bf y})\ldots
\int\limits_{\bf X}\varphi_k({\bf y})
\tilde \nu^{(i_k)}(dt_k,d{\bf y}),
\end{equation}

\vspace{3mm}
\noindent
where permutations $(t_1,\ldots,t_k)$ when summing

$$
\sum_{(t_1,\ldots,t_k)}
$$

\vspace{1mm}
\noindent
in (\ref{jter})
are performed only in the values

\vspace{-1mm}
$$
\varphi_1({\bf y})
\tilde \nu^{(i_1)}(dt_1,d{\bf y})\ldots
\varphi_k({\bf y})
\tilde \nu^{(i_k)}(dt_k,d{\bf y}).
$$

\vspace{2mm}
\noindent
At the same time the indexes near 
upper limits of integration in the iterated stochastic integrals are changed 
correspondently and if $t_r$ swapped with $t_q$ in the  
permutation $(t_1,\ldots,t_k)$, then $i_r$ swapped with $i_q$ in 
the permutation $(i_1,\ldots,i_k)$. Moreover,
$\varphi_r({\bf y})$ swapped with $\varphi_q({\bf y})$
in the permutation $(\varphi_1({\bf y}),\ldots,\varphi_k({\bf y}))$.
Theorem 3 is proved.

Let us consider an example of Theorem 3 usage.
Suppose that $i_1\ne i_2,$ $i_1,i_2=1,\ldots,m$.
According to Theorem 3, we obtain

$$
\int\limits_t^T\int\limits_{\bf X}\varphi_2({\bf y}_2)
\int\limits_{t}^{t_2}\int\limits_{\bf X}\varphi_1({\bf y}_1)
\tilde\nu^{(i_1)}(dt_1,d{\bf y}_1)
\tilde\nu^{(i_2)}(dt_2,d{\bf y}_2)=
$$

\vspace{1mm}
$$
=\frac{T-t}{2}\Biggl(
\pi_{0}^{(1,i_1)}\pi_{0}^{(2,i_2)}
+\sum_{i=1}^{\infty}\frac{1}{\sqrt{4i^2-1}}
\left(\pi_{i-1}^{(1,i_1)}
\pi_{i}^{(2,i_2)}-\pi_{i}^{(1,i_1)}\pi_{i-1}^{(2,i_2)}
\right)\Biggr),
$$

\vspace{3mm}
$$
\int\limits_t^T\int\limits_{\bf X}\varphi_1({\bf y}_1)
\tilde\nu^{(i_1)}(dt_1,d{\bf y}_1)=\sqrt{T-t}\pi_{0}^{(1,i_1)},
$$

\vspace{2mm}
\noindent
where 
$$
\pi_{j}^{(l,i_l)}=
\int\limits_t^T\phi_j(\tau)\int\limits_{\bf X}\varphi_l({\bf y})
\tilde\nu^{(i_l)}(d\tau,d{\bf y})\ \ \ (l=1,\ 2)
$$

\vspace{3mm}
\noindent
and $\{\phi_j(\tau)\}_{j=0}^{\infty}$ is a complete orthonormal
system of Legendre polynomials in the space
$L_2([t, T])$.

\vspace{5mm}

\section{Expansion of Iterated Stochastic Integrals with Respect 
to Martingales}

\vspace{5mm}

Let $(\Omega, {\rm F},{\sf P})$ be a fixed probability
space, let $\{{\rm F}_t,$ $t\in[0,T]\}$
be a non-decreasing family of 
$\sigma$-algebras ${\rm F}_t\subset{\rm F}$, and let
${\rm M}_2(\rho,[0,T])$ be a class 
of ${\rm F}_t$-measurable for each $t\in[0, T]$
martingales $M_t$ satisfying the conditions 

\vspace{-1mm}
\begin{equation}
\label{riss100}
{\sf M}\biggl\{\left(M_s-M_t\right)^2\biggr\}=\int\limits_t^s
\rho(\tau)d\tau,
\end{equation}

$$
{\sf M}\biggl\{\left|M_s-M_t\right|^p\biggr\}\le C_p|s-t|,\ \ \ p=3, 4,\ldots,
$$

\vspace{4mm}
\noindent
where $0\le t< s\le T,$ $\rho(\tau)$ is a non-negative and continuously 
differentiable non-random function at the interval
$[0, T]$, $C_p<\infty$ is a constant.

Let us define the class $H_2(\rho,[0, T])$ 
of stochastic 
processes $\xi_t,$ $t\in[0, T]$ which are 
${\rm F}_t$-measurable for all $t\in[0, T]$ 
and satisfy the condition

\vspace{-1mm}
$$
\int\limits_0^T{\sf M}\left\{\left|\xi_t\right|^2\right\}\rho(t)dt<\infty.
$$

\vspace{2mm}

For any partition 
$\{\tau_j^{(N)}\}_{j=0}^{N}$ of
the interval $[0,T]$ such that

\vspace{2mm}
\begin{equation}
\label{w11}
0=\tau_0^{(N)}<\tau_1^{(N)}<\ldots <\tau_N^{(N)}=T,\ \ \
\max\limits_{0\le j\le N-1}\left|\tau_{j+1}^{(N)}-\tau_j^{(N)}\right|\to\  
0\ \ 
\hbox{if}\ \ N\to \infty
\end{equation}

\vspace{5mm}
\noindent
we will define the sequense of step functions 
$\xi^{(N)}(t,\omega)$ by the following relation

\vspace{2mm}
$$
\xi^{(N)}(t,\omega)=\xi_j\left(\omega\right)\ \ \ 
\hbox{w.~p.~1}\ \ \ \hbox{for}\ \ \ t\in\left[\tau_j^{(N)},
\tau_{j+1}^{(N)}\right),
$$

\vspace{5mm}
\noindent
where $j=0, 1,\ldots,N-1,$\ \  $N=1, 2,\ldots.$

Let us define the stochastic integral with respect to martingale 
from the process $\xi(t,\omega)\in$ $H_2(\rho,[0,T])$ as the following 
mean-square limit \cite{1}

\vspace{1mm}
\begin{equation}
\label{fff}
\hbox{\vtop{\offinterlineskip\halign{
\hfil#\hfil\cr
{\rm l.i.m.}\cr
$\stackrel{}{{}_{N\to \infty}}$\cr
}} }\sum_{j=0}^{N-1}\xi^{(N)}\left(\tau_j^{(N)},\omega\right)
\biggl(M\left(\tau_{j+1}^{(N)},\omega\right)-
M\left(\tau_j^{(N)},\omega\right)\biggr)
\stackrel{\rm def}{=}\int\limits_0^T\xi_\tau dM_\tau,
\end{equation}

\vspace{3mm}
\noindent
where $\xi^{(N)}(t,\omega)$ is any step function
from the class $H_2(\rho,[0,T])$
which converges
to the function $\xi(t,\omega)$
in the following sense

$$
\hbox{\vtop{\offinterlineskip\halign{
\hfil#\hfil\cr
{\rm lim}\cr
$\stackrel{}{{}_{N\to \infty}}$\cr
}} }\int\limits_0^T{\sf M}\left\{\left|\xi^{(N)}(t,\omega)-
\xi(t,\omega)\right|^2\right\}\rho(t)dt=0.
$$

\vspace{3mm}

It is well known  \cite{1} that the stochastic integral

\vspace{-1mm}
$$
\int\limits_0^T\xi_t dM_t
$$

\vspace{2mm}
\noindent
exists and it does not depend on the selection 
of sequence 
$\xi^{(N)}(t,\omega)$ and 
it satisfies w.~p.~1
the following properties 

\vspace{-1mm}
$$
{\sf M}\left\{\int\limits_0^T
\xi_t dM_t\biggl|\biggr.{\rm F}_0\right\}=0,
$$

\vspace{2mm}
$$
{\sf M}\left\{\left|\int\limits_0^T
\xi_t dM_t\right|^2\biggl|\biggr.{\rm F}_0\right\}=
{\sf M}\left\{\int\limits_0^T\xi_t^2\rho(t) dt\biggl|\biggr.{\rm F}_0\right\},
$$

\vspace{1mm}
$$
\int\limits_0^T(\alpha\xi_t+\beta\psi_t)dM_t=
\alpha\int\limits_0^T\xi_t dM_t+\beta
\int\limits_0^T\psi_t dM_t,
$$

\vspace{4mm}
\noindent
where $\xi_t,$ $\phi_t\in H_2(\rho,[0, T]),$\ \
$\alpha,\ \beta\in\mathbb{R}^1.$

Let $Q_4(\rho,[0,T])$ be the class 
of martingales $M_t,$ $t\in[0,T]$ 
which satisfy the following conditions:

1. $M_t,$ $t\in[0,T]$ belongs to the class
${\rm M}_2(\rho,[0,T]).$ 

2. For some 
$\alpha>0$ the following estimate is correct

\begin{equation}
\label{durra1}
{\sf M}\left\{\left|\int\limits_{t}^{\tau}
g(s) dM_s\right|^4\right\}\le K_4 \int\limits_{t}^{\tau}
|g(s)|^{\alpha}ds,
\end{equation}

\vspace{2mm}
\noindent
where $0\le t<\tau\le T,$\ \
$g(s)$ is a bounded non-random function 
at the interval $[0, T]$,\ \
$K_4<\infty$ is a constant.

Let $G_n(\rho,[0,T])$ be the class 
of martingales $M_t,$ $t\in[0,T]$ 
which satisfy the following conditions:

1. $M_t,$ $t\in[0,T]$ belongs to the class
${\rm M}_2(\rho,[0,T]).$ 

2. The following estimate is correct

$$
{\sf M}\left\{\left|\int\limits_{t}^{\tau}
g(s) dM_s\right|^n\right\} < \infty,
$$

\vspace{3mm}
\noindent
where $0\le t<\tau\le T,$\ \ $n\in \mathbb{N},$\ \ $g(s)$ 
is the same function as in the definition of
$Q_4(\rho,[0,T])$.

Let us 
remind that if $\left(\xi_t\right)^n\in H_2(\rho,[0,T])$
with $\rho(t)\equiv 1,$ then the following estimate is correct
\cite{1}

\vspace{-1mm}
\begin{equation}
\label{1.5aa}
{\sf M}\left\{\left|\int\limits_{t}^{\tau}
\xi_s ds\right|^{2n}\right\}\le (\tau-t)^{2n-1}
\int\limits_{t}^{\tau}
{\sf M}\left\{\left|\xi_s\right|^{2n}\right\}ds,\ \ \ 0\le t<\tau\le T.
\end{equation}

\vspace{2mm}

Let us consider the iterated stochastic integral with respect to martingales

\vspace{-1mm}
\begin{equation}
\label{mart}
J[\psi^{(k)}]_{T,t}^M=
\int\limits_t^T\psi_k(t_k)\ldots
\int\limits_t^{t_2}\psi_1(t_1)
dM_{t_1}^{(1,i_1)}\ldots dM_{t_k}^{(k,i_k)}\ \ \ 
(i_1,\ldots,i_k=0, 1,\ldots,m),
\end{equation}

\vspace{2mm}
\noindent
where every $\psi_l(\tau)$ $(l=1,\ldots, k)$ is a continuous 
non-random function
at the interval $[t, T],$\ 
$M^{(r,i)}$ $(r=1,\ldots,k)$ are independent martingales 
for various $i=1,\ldots,m,$\ 
$M_{\tau}^{(r,0)}\stackrel{\rm def}{=}\tau.$

Let us formulate the following theorem. 

\vspace{2mm}

{\bf Theorem 4}\ \cite{2017}-\cite{12aa-afterxxx}.\ 
{\it Suppose that the following 
conditions are fulfilled{\rm :}

{\rm 1}.\ Every $\psi_l(\tau)\ (l=1,\ldots,k)$ is a 
continuous non-random function at
the interval $[t, T]$.

{\rm 2}.\ $\{\phi_j(x)\}_{j=0}^{\infty}$ is a complete orthonormal
system of functions in the space 
$L_2([t,T]),$ 
each function of which for finite $j$ satisfies the condition 
$(\star)$ {\rm(}see Sect. {\rm 1)}.

{\rm 3}. $M_{\tau}^{(l,i_l)}\in Q_4(\rho,[t,T]),$ $G_n(\rho,[t,T])$
with
$n=2^{k+1},$ 
$i_l=1,\ldots,m,$\ \  $l=1,\ldots,k.$

Then, for the iterated stochastic integral 
$J[\psi^{(k)}]_{T,t}^M$ with respect to martingales 
defined by {\rm (\ref{mart})}
the following expansion

\vspace{1mm}
$$
J[\psi^{(k)}]_{T,t}^M=
\hbox{\vtop{\offinterlineskip\halign{
\hfil#\hfil\cr
{\rm l.i.m.}\cr
$\stackrel{}{{}_{p_1,\ldots,p_k\to \infty}}$\cr
}} }\sum_{j_1=0}^{p_1}\ldots\sum_{j_k=0}^{p_k}
C_{j_k\ldots j_1}\Biggl(
\prod_{l=1}^k\xi_{j_l}^{(l,i_l)}-
\Biggr.
$$

\vspace{3mm}
$$
-\Biggl.
\hbox{\vtop{\offinterlineskip\halign{
\hfil#\hfil\cr
{\rm l.i.m.}\cr
$\stackrel{}{{}_{N\to \infty}}$\cr
}} }\sum_{(l_1,\ldots,l_k)\in {\rm G}_k}
\phi_{j_{1}}(\tau_{l_1})
\Delta{M}_{\tau_{l_1}}^{(1,i_1)}\ldots
\phi_{j_{k}}(\tau_{l_k})
\Delta{M}_{\tau_{l_k}}^{(k,i_k)}\Biggr)
$$

\vspace{6mm}
\noindent
converging in the mean-square sense is valid, where $i_1,\ldots,i_k=0,1,\ldots,m,$\ 
$\left\{\tau_{j}\right\}_{j=0}^{N}$ is a partition of
the interval $[t, T]$ which satisfies the condition {\rm (\ref{1111}),}
$\Delta{M}_{\tau_{j}}^{(r,i)}=
M_{\tau_{j+1}}^{(r,i)}-M_{\tau_{j}}^{(r,i)}$
$(i=0, 1,\ldots,m,\ \ r=1,\ldots,k),$

$$
{\rm G}_k={\rm H}_k\backslash{\rm L}_k,\ \ \
{\rm H}_k=\{(l_1,\ldots,l_k):\ l_1,\ldots,l_k=0,\ 1,\ldots,N-1\},
$$

\vspace{-1mm}
$$
{\rm L}_k=\{(l_1,\ldots,l_k):\ l_1,\ldots,l_k=0,\ 1,\ldots,N-1;\
l_g\ne l_r\ (g\ne r);\ g, r=1,\ldots,k\},
$$

\vspace{5mm}
\noindent
${\rm l.i.m.}$ is a limit in the mean-square sense$,$\

$$
\xi_{j}^{(l,i_l)}=
\int\limits_t^T \phi_{j}(s) d{M}_s^{(l,i_l)}
$$

\vspace{2mm}
\noindent
are independent for various
$i_l=1,\ldots,m,$\ $l=1,\ldots,k$
and uncorrelated for various $j$
{\rm (}if $\rho(\tau)$ is a constant, $i_l\ne 0${\rm )} random variables,

\vspace{1mm}
$$
C_{j_k\ldots j_1}=\int\limits_{[t,T]^k}
K(t_1,\ldots,t_k)\prod_{l=1}^{k}\phi_{j_l}(t_l)dt_1\ldots dt_k
$$

\vspace{2mm}
\noindent
is the Fourier coefficient,

$$
K(t_1,\ldots,t_k)=
\begin{cases}
\psi_1(t_1)\ldots \psi_k(t_k)\ &\hbox{for}\ \ t_1<\ldots<t_k\\
~\\
~\\
0\ &\hbox{otherwise}
\end{cases},\ \ \ \ t_1,\ldots,t_k\in[t, T],\ \ \ \ k\ge 2
$$

\vspace{5mm}
\noindent
and 
$K(t_1)\equiv\psi_1(t_1)$ for $t_1\in[t, T].$  
}

\vspace{4mm}

{\bf Remark 1.}\ {\it Note that from Theorem {\rm 4} 
for the case $\rho(\tau)\equiv 1$ we obtain the variant of Theorem 
{\rm 1}.}

\vspace{2mm}

{\bf Proof.} The scheme of the proof of Theorem 4 is the same with the scheme
of the proof of Theorem 1 in \cite{arxiv-1}
or Theorem 1.1 in
\cite{2018a}-\cite{12aa-afterxxx}
(also see \cite{2006}-\cite{2013}, \cite{arxiv-1}).
Some differences will take place in 
the proof of the following lemmas (Lemmas 3, 4) and in the final part 
of the proof of Theorem 4.

\vspace{2mm}

{\bf Lemma 3.} {\it Suppose that
$M_{\tau}^{(r,i)}\in{\rm M}_2(\rho,[t,T]),$\ $M_{\tau}^{(r,0)}=\tau$\
$(i=0, 1,\ldots,m,$\ \  $r=1,\ldots,k),$ and
every $\psi_l(\tau)$ $(l=1,\ldots, k)$ is a continuous non-random
function at the interval $[t,T]$.
Then

\vspace{-2mm}
\begin{equation}
\label{1.9aa}
J[\psi^{(k)}]_{T,t}^M=
\hbox{\vtop{\offinterlineskip\halign{
\hfil#\hfil\cr
{\rm l.i.m.}\cr
$\stackrel{}{{}_{N\to \infty}}$\cr
}} }
\sum_{j_k=0}^{N-1}
\ldots \sum_{j_1=0}^{j_{2}-1}
\prod_{l=1}^k \psi_{l}(\tau_{j_l})\Delta M_{\tau_{j_l}}^{(l,i_l)}\ \ \
\hbox{w. p. {\rm 1}},
\end{equation}

\vspace{3mm}
\noindent
where 
$\{\tau_j\}_{j=0}^N$ is a partition of the interval
$[0,T]$ satisfying the condition {\rm(\ref{1111}).}
}

\vspace{2mm}

{\bf Proof.} According to properties of the stochastic integral with respect
to martingale, we have \cite{1}

\vspace{-1mm}
\begin{equation}
\label{u1}
{\sf M}\left\{\left(\int\limits_{t}^{\tau}
\xi_s dM_s^{(l,i_l)}\right)^2\right\}=
\int\limits_{t}^{\tau}{\sf M}\left\{\left|\xi_s\right|^2\right\}\rho(s)ds,
\end{equation}

\begin{equation}
\label{u2}
{\sf M}\left\{\left(\int\limits_{t}^{\tau}
\xi_s ds\right)^2\right\}\le (\tau-t)
\int\limits_{t}^{\tau}{\sf M}\left\{\left|\xi_s\right|^2\right\}ds,
\end{equation}

\vspace{2mm}
\noindent
where 
$\xi_s\in H_2(\rho,[0,T]),$\ $0\le t<\tau\le T,$\ \
$i_l=1,\ldots,m,$\ \  $l=1,\ldots,k$. 
Then the integral sum 
of the integral 
$J[\psi^{(k)}]_{T,t}^M$
under
the conditions of Lemma 3 can be represented as a sum of the 
expression from the right-hand side of (\ref{1.9aa}) before passing to the 
limit and the 
value which 
converges
to zero 
in the mean-square sense if 
$N\to \infty.$ 
More detailed proof of the analogous lemma for the case $\rho(\tau)\equiv 1$ 
can be found in \cite{2006}-\cite{2013}, \cite{arxiv-1}.

In the case when the functions 
$\psi_l(\tau)$ $(l=1,\ldots,k)$ satisfy the condition
$(\star)$ {\rm (}see Sect. {\rm 1)}
we can suppose that among the points
$\tau_j,$ $j=0,1,\ldots,N$ there are all points of 
jumps of the functions $\psi_l(\tau)$ 
$(l=1,\ldots,k).$ 
Then can apply the argumentation as in Sect.~4 from
\cite{arxiv-1} (also see \cite{2006}-\cite{12aa-afterxxx}).

Let us define the folloing multiple stochastic integral

\begin{equation}
\label{777666}
\hbox{\vtop{\offinterlineskip\halign{
\hfil#\hfil\cr
{\rm l.i.m.}\cr
$\stackrel{}{{}_{N\to \infty}}$\cr
}} }
\sum_{j_1,\ldots,j_k=0}^{N-1}
\Phi(\tau_{j_1},\ldots,\tau_{j_k})
\prod_{l=1}^k 
\Delta M_{\tau_{j_l}}^{(l,i_l)}
\stackrel{\rm def}{=}I[\Phi]_{T,t}^{(k)},
\end{equation}

\vspace{3mm}
\noindent
where $\{\tau_j\}_{j=0}^N$ is a partition of the interval 
$[0,T]$ satisfying the condition (\ref{1111})
and
$\Phi(t_1,\ldots,t_k):$ $[t, T]^k\to\mathbb{R}^1$ is a bounded non-random
function.

\vspace{2mm}

{\bf Lemma 4.} {\it Suppose that 
$M_s^{(l,i_l)}\in Q_4(\rho,[t, T]),$\ $G_n(\rho,[t,T])$\ with
$n=2^{k+1},$\ $k\in\mathbb{N}$\   $(i_l=0, 1,\ldots,m,$\ \ 
$l=1,\ldots,k)$
and the functions $g_1(s),\ldots, g_k(s)$ satisfy the condition 
$(\star)$ {\rm(}see Sect. {\rm 1)}.
Then

$$
\prod_{l=1}^k \int\limits_t^T g_l(s) 
dM_s^{(l,i_l)}=
I[\Phi]_{T,t}^{(k)}\ \ \ \hbox{w.~p.~{\rm 1}},
$$

\vspace{2mm}
\noindent
where
$$
\Phi(t_1,\ldots,t_k)=\prod\limits_{l=1}^k g_l(t_l).
$$
}

\vspace{2mm}

{\bf Proof.} Let us denote

$$
J[g_l]_N\stackrel{\rm def}{=}\sum\limits_{j=0}^{N-1}
g_l(\tau_j)\Delta M_{\tau_j}^{(l,i_l)},\ \ \
J[g_l]_{T,t}\stackrel{\rm def}{=}\int\limits_t^T g_l(s)
dM_s^{(l,i_l)},
$$

\vspace{4mm}
\noindent
where $\{\tau_j\}_{j=0}^N$ is a partition of the interval 
$[t,T]$ satisfying the condition (\ref{1111}).

Note that

$$
\prod_{l=1}^k J[g_l]_N-\prod_{l=1}^k J[g_l]_{T,t}
=
$$

\vspace{2mm}
$$
=\sum_{l=1}^k \left(\prod_{q=1}^{l-1} J[g_q]_{T,t}\right)
\left(J[g_l]_N-
J[g_l]_{T,t}\right)\left(\prod_{q=l+1}^k J[g_q]_N\right).
$$

\vspace{5mm}

Using the Minkowski inequality and the
inequality of Cauchy-Bu\-nya\-kov\-sky
as well as the conditions of Lemma 4, we obtain

\begin{equation}
\label{2000.4.300}
\left({\sf M}\left\{\left|\prod_{l=1}^k J[g_l]_N-
\prod_{l=1}^k J[g_l]_{T,t}\right|^2\right\}
\right)^{1/2}\le C_k
\sum_{l=1}^k
\left({\sf M}
\left\{\biggl|J[g_l]_N-J[g_l]_{T,t}
\biggr|^4\right\}\right)^{1/4},
\end{equation}

\vspace{4mm}
\noindent
where $C_k<\infty$ is a constant.

We have

\vspace{-1mm}
$$
J[g_l]_N-J[g_l]_{T,t}
=\sum\limits_{q=0}^{N-1}J[\Delta g_{l}]_{\tau_{q+1},\tau_q},
$$

$$
J[\Delta g_{l}]_{\tau_{q+1},\tau_q}
=\int\limits_{\tau_q}^{\tau_{q+1}}\left(
g_l(\tau_q)-g_l(s)\right)
dM_s^{(l,i_l)}.
$$

\vspace{4mm}

Let as introduce the notation

$$
g_l^{(N)}(s)=g_l(\tau_q),\ \ \ s\in [\tau_q, \tau_{q+1}), 
\ \ \ q=0, 1, \ldots,N-1.
$$

\vspace{3mm}

Then
$$
J[\bar g_l]_N-J[\bar g_l]_{T,t}
=\sum\limits_{q=0}^{N-1}J[\Delta\bar g_{l}]_{\tau_{q+1},\tau_q}=
$$

$$
=
\int\limits_{t}^{T}\left(g_l^{(N)}(s)-
g_l(s)\right)dM_s^{(l,i_l)}.
$$

\vspace{3mm}

Applying the estimate (\ref{durra1}), 
we obtain

$$
{\sf M}\left\{\left|
\int\limits_{t}^{T}\left(g_l^{(N)}(s)-
g_l(s)\right)dM_s^{(l,i_l)}\right|^4\right\}\le
K_4\int\limits_t^T\left|
g_l^{(N)}(s)-
g_l(s)\right|^{\alpha}ds=
$$
$$
=K_4\sum_{q=0}^{N-1}\int\limits_{\tau_q}^{\tau_{q+1}}\left|
g_l(\tau_q)-
g_l(s)\right|^{\alpha}ds< K_4 \varepsilon^{\alpha}\ 
\sum_{q=0}^{N-1}\left(\tau_{q+1}-\tau_q\right)=
$$
\begin{equation}
\label{durra2}
=K_4 \varepsilon^{\alpha}(T-t).
\end{equation}

\vspace{3mm}

Note that deriving (\ref{durra2}) 
we used 
the estimate

\begin{equation}
\label{durra3}
\left|g_l(\tau_q)-
g_l(s)\right|<\varepsilon,\ \ \ s\in[\tau_q,\tau_{q+1}],\ \ \ q=0,1,\ldots,N-1,
\end{equation}

\vspace{3mm}
\noindent
where $\varepsilon$ is an arbitrary small positive real number.

Note that (\ref{durra3}) is valid
if the functions $g_l(s)$ are continuous at the interval
$[t, T]$, i.e. these functions are uniformly continuous at this interval.
So, $\left|g_l(\tau_q)-g_l(s)\right|<\varepsilon$
if $s\in [\tau_q, \tau_{q+1}],$ where
$|\tau_{q+1}-\tau_q|<\delta(\varepsilon),$ $q=0, 1,\ldots,N-1$
($\delta(\varepsilon)>0$ exists
for any $\varepsilon>0$ and it does not
depend on points of the interval $[t, T]$).

Thus, taking into account (\ref{durra2}),
we obtain that the right-hand side of (\ref{2000.4.300}) 
converges
to zero when $N\to\infty.$
Considering this fact, we come to 
the statement of Lemma 4.

In the case when the functions 
$g_l(s)$ 
$(l=1,\ldots,k)$ satisfy the condition
$(\star)$ {\rm (}see Sect.~{\rm 1)}
we can suppose that among the points
$\tau_q,$ $q=0,1,\ldots,N$ there are all points of 
jumps of the functions $g_l(s)$ 
$(l=1,\ldots,k)$. Further, 
we can apply the argumentation as in Sect.~4 from
\cite{arxiv-1} (also see \cite{2006}-\cite{12aa-afterxxx}).

Obviously, if $i_l=0$ for some $l=1,\ldots,k,$ then
we also come to the 
statement of Lemma 4 with using (\ref{u2}).
Lemma 4 is proved.

Proving Theorem 4 according to the scheme used for the proof
of Theorem 1 in \cite{arxiv-1} 
or Theorem 1.1 in
\cite{2018a}-\cite{12aa-afterxxx} (also see \cite{2006} (Theorem 5.1, P. 236-237), 
\cite{2017} (Theorem 1, P. A.22-A.23),
\cite{2017-1} (Theorem 5.1, P. A.250), 
\cite{2018} (Theorem 5.1, P. A.252-A.253)) and
using Lemmas 3, 4 together with the estimates (\ref{u1}), (\ref{u2})
for moments
of stochastic integrals with respect to martingales, we obtain

\vspace{1mm}
$$
{\sf M}\left\{\left(R_{T,t}^{p_1,\ldots,p_k}\right)^2\right\}
\le 
$$

\vspace{2mm}
$$
\le
C_k
\sum_{(t_1,\ldots,t_k)}
\int\limits_{t}^{T}
\ldots
\int\limits_{t}^{t_2}
\left(K(t_1,\ldots,t_k)-
\sum_{j_1=0}^{p_1}\ldots
\sum_{j_k=0}^{p_k}
C_{j_k\ldots j_1}
\prod_{l=1}^k\phi_{j_l}(t_l)\right)^2\times
$$

\vspace{2mm}
\begin{equation}
\label{z2}
\times{\tilde \rho}_1(t_1)dt_1
\ldots
{\tilde \rho}_k(t_k)dt_k\le
\end{equation}

\vspace{3mm}
$$
\le
\bar{C_k}
\sum_{(t_1,\ldots,t_k)}
\int\limits_{t}^{T}
\ldots
\int\limits_{t}^{t_2}
\left(K(t_1,\ldots,t_k)-
\sum_{j_1=0}^{p_1}\ldots
\sum_{j_k=0}^{p_k}
C_{j_k\ldots j_1}
\prod_{l=1}^k\phi_{j_l}(t_l)\right)^2dt_1\ldots dt_k=
$$

\vspace{5mm}
$$
=\bar{C_k}
\int\limits_{[t,T]^k}
\Biggl(K(t_1,\ldots,t_k)-
\sum_{j_1=0}^{p_1}\ldots
\sum_{j_k=0}^{p_k}
C_{j_k\ldots j_1}
\prod_{l=1}^k\phi_{j_l}(t_l)\Biggr)^2
dt_1
\ldots
dt_k \to 0
$$

\vspace{5mm}
\noindent
when $p_1,\ldots,p_k\to\infty,$
where constant $\bar{C_k}$ depends only on $k$ (multiplicity of 
the iterated stochastic integral with respect to martingales) and  
${\tilde \rho}_l(s)\equiv\rho(s)$ or ${\tilde \rho}_l(s)\equiv 1$
$(l=1,\ldots,k)$. 
At that permutations $(t_1,\ldots,t_k)$ when summing

\vspace{-1mm}
$$
\sum_{(t_1,\ldots,t_k)}
$$

\vspace{2mm}
\noindent
in (\ref{z2})
are performed only in the values $dt_1\ldots dt_k$ and
indexes near 
upper limits of integration are changed 
correspondently.
Moreover, $R_{T,t}^{p_1,\ldots,p_k}$ has the following 
form

\vspace{1mm}
$$
R_{T,t}^{p_1,\ldots,p_k}
=\sum_{(t_1,\ldots,t_k)}
\int\limits_{t}^{T}
\ldots
\int\limits_{t}^{t_2}
\left(K(t_1,\ldots,t_k)-
\sum_{j_1=0}^{p_1}\ldots
\sum_{j_k=0}^{p_k}
C_{j_k\ldots j_1}
\prod_{l=1}^k\phi_{j_l}(t_l)\right)\times
$$

\vspace{1mm}
\begin{equation}
\label{jter1}
\times
dM_{t_1}^{(1,i_1)}\ldots
dM_{t_k}^{(k,i_k)},
\end{equation}

\vspace{5mm}
\noindent
where permutations $(t_1,\ldots,t_k)$ when summing

\vspace{-1mm}
$$
\sum_{(t_1,\ldots,t_k)}
$$

\vspace{2mm}
\noindent
in (\ref{jter1})
are performed only in the values

$$
dM_{t_1}^{(1,i_1)}\ldots
dM_{t_k}^{(k,i_k)}.
$$

\vspace{3mm}
\noindent
At the same time the indexes near 
upper limits of integration in the iterated stochastic integrals are changed 
correspondently and if $t_r$ swapped with $t_q$ in the  
permutation $(t_1,\ldots,t_k)$, then $i_r$ swapped with $i_q$ in 
the permutation $(i_1,\ldots,i_k)$. Moreover,
$r$ swapped with $q$
in the permutation $(1,\ldots,k)$.
Theorem 4 is proved.

\vspace{5mm}

\section{Expansion of Iterated Ito Stochastic Integrals 
Based on Generalized Multiple Fourier Series.
The Case of Complete Orthonormal With We\-ight $r(t_1)\ldots r(t_k)\ge 0$  
Systems of Functions 
in the Space $L_2([t, T]^k)$}

\vspace{5mm}

In this section, we consider modifications of Theorems 1, 2 for 
the case of complete orthonormal with weight $r(t_1)\ldots r(t_k)\ge 0$ 
systems of functions 
in the space $L_2([t, T]^k),$ $k\in\mathbb{N}$.

Let $\{\Psi_j(x)\}_{j=0}^{\infty}$ be a complete orthonormal 
with weight $r(x)\ge 0$ 
system of functions in the space $L_2([t, T]).$ It is well known that the
Fourier
series with respect to the system 

\vspace{-1mm}
$$
\{\Psi_j(x)\}_{j=0}^{\infty}
$$

\vspace{3mm}
\noindent
of the function $f(x)$ $\left(f(x)\sqrt{r(x)}\in L_2([t, T])\right)$ 
converges 
to the function $f(x)$ in the
mean-square sense with weight $r(x),$ i.e.

\vspace{-1mm}
\begin{equation}
\label{g1}
\lim\limits_{p\to\infty}
\int\limits_t^T\biggl(f(x)-\sum\limits_{j=0}^p 
{\tilde C}_j \Psi_j(x)\biggr)^2 r(x)dx = 0,
\end{equation}

\vspace{3mm}
\noindent
where
\begin{equation}
\label{h1}
{\tilde C}_j=\int\limits_t^T f(x)\Psi_j(x)r(x)dx
\end{equation}

\vspace{2mm}
\noindent
is the Fourier coefficient.

Obviously, the relation (\ref{g1}) can be obtained if we will 
expand the function
$f(x)\sqrt{r(x)}\in L_2([t, T])$ into a usual Fourier series with respect
to the complete orthonormal with weight $1$ system of functions
$$
\left\{\Psi_j(x)\sqrt{r(x)}\right\}_{j=0}^{\infty}
$$ 

\noindent
in 
the space $L_2([t, T]).$ Then
$$
\lim\limits_{p\to\infty}
\int\limits_t^T\biggl(f(x)\sqrt{r(x)}-\sum\limits_{j=0}^p {\tilde C}_j 
\Psi_j(x)\sqrt{r(x)}\biggr)^2dx = 
$$
\begin{equation}
\label{g2}
=\lim\limits_{p\to\infty}
\int\limits_t^T\biggl(f(x)-
\sum\limits_{j=0}^p {\tilde C}_j \Psi_j(x)\biggr)^2 r(x)dx = 0,
\end{equation}

\noindent
where ${\tilde C}_j$ has the form (\ref{h1}).

Let us consider an obvious generalization of this approach 
to the case of several
variables.
Let us expand the function $K(t_1,\ldots,t_k)$ such that

$$
K(t_1,\ldots,t_k)\prod\limits_{l=1}^k \sqrt{r(t_l)}\in L_2([t, T]^k)
$$

\vspace{3mm}
\noindent
using the complete orthonormal system of functions 

\vspace{-1mm}
$$
\prod\limits_{l=1}^k \Psi_{j_l}(t_l)\sqrt{r(t_l)},\ \ \ 
j_l=0, 1, 2, \ldots,\ \ \  l=1,\ldots,k
$$

\vspace{3mm}
\noindent
in the space $L_2([t, T]^k)$ into the generalized multiple Fourier 
series. 

It is well known that the mentioned
generalized multiple Fourier series converges in the mean-square sense,
i.e.

$$
\lim\limits_{p_1,\ldots,p_k\to\infty}
\int\limits_{[t,T]^k}
\left(K(t_1,\ldots,t_k)\prod\limits_{l=1}^k \sqrt{r(t_l)}-
\sum\limits_{j_1=0}^{p_1}\ldots\sum\limits_{j_k=0}^{p_k}
{\tilde C}_{j_k\ldots j_1}
\prod\limits_{l=1}^k \Psi_{j_l}(t_l)\sqrt{r(t_l)}\right)^2\times
$$

$$
\times
dt_1\ldots dt_k=
$$

\vspace{1mm}

\begin{equation}
\label{z1aaa}
=\lim\limits_{p_1,\ldots,p_k\to\infty}
\int\limits_{[t,T]^k}
\left(K(t_1,\ldots,t_k)-
\sum\limits_{j_1=0}^{p_1}\ldots\sum\limits_{j_k=0}^{p_k}
{\tilde C}_{j_k\ldots j_1}\prod\limits_{l=1}^k \Psi_{j_l}(t_l)\right)^2 
\times
$$

\vspace{1mm}
$$
\times
\left(\prod\limits_{l=1}^k r(t_l)\right)
dt_1\ldots dt_k=0,
\end{equation}

\vspace{5mm}
\noindent
where
$$
{\tilde C}_{j_k\ldots j_1}=\int\limits_{[t,T]^k}
K(t_1,\ldots,t_k)\prod\limits_{l=1}^k 
\biggl(\Psi_{j_l}(t_l)r(t_l)\biggr)dt_1\ldots dt_k.
$$

\vspace{3mm}

Let us consider 
the following iterated Ito 
stochastic integrals

\vspace{-1mm}
\begin{equation}
\label{ito-rr}
{\tilde J}[\psi^{(k)}]_{T,t}=\int\limits_t^T\psi_k(t_k)\sqrt{r(t_k)} 
\ldots \int\limits_t^{t_{2}}
\psi_1(t_1)\sqrt{r(t_1)} d{\bf w}_{t_1}^{(i_1)}\ldots
d{\bf w}_{t_k}^{(i_k)},
\end{equation}

\vspace{2mm}
\noindent
where every $\psi_l(\tau)$ $(l=1,\ldots,k)$ is 
a non-random function on $[t, T]$,
${\bf w}_{\tau}^{(i)}={\bf f}_{\tau}^{(i)}$
for $i=1,\ldots,m,$
${\bf w}_{\tau}^{(0)}=\tau,$ and
$i_1,\ldots,i_k=0, 1,\ldots,m.$

So, we obtain the following modification of Theorem 1.

\vspace{2mm}

{\bf Theorem 5}\ \cite{2018}-\cite{2018aaa}, \cite{arxiv-26b}.\
{\it Suppose that
every $\psi_l(\tau)$ $(l=$ $1,\ldots, k)$ is a continuous 
non-random function on 
$[t, T]$ and
$\{\Psi_j(x)\sqrt{r(x)}\}_{j=0}^{\infty}$ $(r(x)\ge 0)$
is a complete orthonormal 
system of functions in the space $L_2([t,T]),$ each function 
$\Psi_j(x)\sqrt{r(x)}$
of which 
for finite $j$ satisfies the condition 
$(\star)$ {\rm (}see Sect.~{\rm 1)}.
Then

$$
{\tilde J}[\psi^{(k)}]_{T,t} =
\hbox{\vtop{\offinterlineskip\halign{
\hfil#\hfil\cr
{\rm l.i.m.}\cr
$\stackrel{}{{}_{p_1,\ldots,p_k\to \infty}}$\cr
}} }\sum_{j_1=0}^{p_1}\ldots\sum_{j_k=0}^{p_k}
{\tilde C}_{j_k\ldots j_1}\Biggl(
\prod_{l=1}^k{\tilde \zeta}_{j_l}^{(i_l)} -
\Biggr.
$$

\vspace{2mm}
\begin{equation}
\label{tyyy-rr}
-\Biggl.
\hbox{\vtop{\offinterlineskip\halign{
\hfil#\hfil\cr
{\rm l.i.m.}\cr
$\stackrel{}{{}_{N\to \infty}}$\cr
}} }\sum_{(l_1,\ldots,l_k)\in {\rm G}_k}
\Psi_{j_{1}}(\tau_{l_1})\sqrt{r(\tau_{l_1})}
\Delta{\bf w}_{\tau_{l_1}}^{(i_1)}\ldots
\Psi_{j_{k}}(\tau_{l_k})\sqrt{r(\tau_{l_k})}
\Delta{\bf w}_{\tau_{l_k}}^{(i_k)}\Biggr),
\end{equation}

\vspace{7mm}
\noindent
where
$$
{\rm G}_k={\rm H}_k\backslash{\rm L}_k,\ \ \
{\rm H}_k=\left\{(l_1,\ldots,l_k):\ l_1,\ldots,l_k=0,\ 1,\ldots,N-1\right\},
$$

\vspace{-2mm}
$$
{\rm L}_k=\left\{(l_1,\ldots,l_k):\ l_1,\ldots,l_k=0,\ 1,\ldots,N-1;\
l_g\ne l_r\ (g\ne r);\ g, r=1,\ldots,k\right\},
$$

\vspace{3mm}
\noindent
${\rm l.i.m.}$ is a limit in the mean-square sense,
$i_1,\ldots,i_k=0,1,\ldots,m,$ 
$$
{\tilde \zeta}_{j}^{(i)}=
\int\limits_t^T \Psi_{j}(s)\sqrt{r(s)}d{\bf w}_s^{(i)}
$$

\vspace{2mm}
\noindent
are independent standard Gaussian random variables
for various
$i$ or $j$ {\rm(}in the case when $i\ne 0${\rm),}
$\Delta{\bf w}_{\tau_{j}}^{(i)}=
{\bf w}_{\tau_{j+1}}^{(i)}-{\bf w}_{\tau_{j}}^{(i)}$
$(i=0, 1,\ldots,m),$
$\left\{\tau_{j}\right\}_{j=0}^{N}$ is a partition of
$[t,T]$ which satisfies the condition {\rm (\ref{1111})},
\begin{equation}
\label{koef}
{\tilde C}_{j_k\ldots j_1}=\int\limits_{[t,T]^k}
K(t_1,\ldots,t_k)
\prod_{l=1}^{k}\biggl(\Psi_{j_l}(t_l)r(t_l)\biggr)dt_1\ldots dt_k
\end{equation}

\vspace{2mm}
\noindent
is the Fourier coefficient,

$$
K(t_1,\ldots,t_k)=
\begin{cases}
\psi_1(t_1)\ldots \psi_k(t_k)\ &\hbox{for}\ \ t_1<\ldots<t_k\\
~\\
~\\
0\ &\hbox{otherwise}
\end{cases},\ \ \ \ t_1,\ldots,t_k\in[t, T],\ \ \ \ k\ge 2
$$

\vspace{5mm}
\noindent
and 
$K(t_1)\equiv\psi_1(t_1)$ for $t_1\in[t, T].$  
}

\vspace{2mm}

{\bf Proof.}
According to Lemmas 1--3 in \cite{arxiv-1} or
Lemmas 1.1--1.3 in \cite{2018a}-\cite{12aa-afterxxx}
(also see \cite{2006}-\cite{2018}),
we get the following representation w.~p.~1 

$$
{\tilde J}[\psi^{(k)}]_{T,t}=
$$

\vspace{2mm}
$$
=\sum_{(t_1,\ldots,t_k)}
\int\limits_{t}^{T}
\ldots
\int\limits_{t}^{t_2}
K(t_1,\ldots,t_k)\prod\limits_{l=1}^k \sqrt{r(t_l)}d{\bf w}_{t_1}^{(i_1)}
\ldots
d{\bf w}_{t_k}^{(i_k)}=
$$

\vspace{6mm}
$$
=
\sum_{j_1=0}^{p_1}\ldots
\sum_{j_k=0}^{p_k}
{\tilde C}_{j_k\ldots j_1}
\sum_{(t_1,\ldots,t_k)}\int\limits_{t}^{T}
\ldots
\int\limits_{t}^{t_2}
\prod\limits_{l=1}^k\left(\Psi_{j_l}(t_l)\sqrt{r(t_l)}\right)
d{\bf w}_{t_1}^{(i_1)}
\ldots
d{\bf w}_{t_k}^{(i_k)}
+
$$

\vspace{2mm}
$$
+{\tilde R}_{T,t}^{p_1,\ldots,p_k}=
$$

\vspace{6mm}
$$
=\sum_{j_1=0}^{p_1}\ldots
\sum_{j_k=0}^{p_k}
{\tilde C}_{j_k\ldots j_1}\times
$$

\vspace{2mm}
$$
\times             
\hbox{\vtop{\offinterlineskip\halign{
\hfil#\hfil\cr
{\rm l.i.m.}\cr
$\stackrel{}{{}_{N\to \infty}}$\cr
}} }
\sum\limits_{\stackrel{l_1,\ldots,l_k=0}{{}_{l_q\ne l_r;\ 
q\ne r;\ q, r=1,\ldots, k}}}^{N-1}
\Psi_{j_1}(\tau_{l_1})\sqrt{r(\tau_{l_1})}\Delta{\bf w}_{\tau_{l_1}}^{(i_1)}
\ldots
\Psi_{j_k}(\tau_{l_k})\sqrt{(\tau_{l_k})}
\Delta{\bf w}_{\tau_{l_k}}^{(i_k)}+
$$

\vspace{2mm}
$$
+{\tilde R}_{T,t}^{p_1,\ldots,p_k}=
$$

\vspace{6mm}
$$
=\sum_{j_1=0}^{p_1}\ldots
\sum_{j_k=0}^{p_k}
{\tilde C}_{j_k\ldots j_1}\times
$$

\vspace{2mm}
$$
\times\left(
\hbox{\vtop{\offinterlineskip\halign{
\hfil#\hfil\cr
{\rm l.i.m.}\cr
$\stackrel{}{{}_{N\to \infty}}$\cr
}} }\sum_{l_1,\ldots,l_k=0}^{N-1}
\Psi_{j_1}(\tau_{l_1})\sqrt{r(\tau_{l_1})}\Delta{\bf w}_{\tau_{l_1}}^{(i_1)}
\ldots
\Psi_{j_k}(\tau_{l_k})\sqrt{(\tau_{l_k})}
\Delta{\bf w}_{\tau_{l_k}}^{(i_k)}
-\right.
$$

\vspace{2mm}
$$
-\left.
\hbox{\vtop{\offinterlineskip\halign{
\hfil#\hfil\cr
{\rm l.i.m.}\cr
$\stackrel{}{{}_{N\to \infty}}$\cr
}} }\sum_{(l_1,\ldots,l_k)\in {\rm G}_k}
\Psi_{j_1}(\tau_{l_1})\sqrt{r(\tau_{l_1})}\Delta{\bf w}_{\tau_{l_1}}^{(i_1)}
\ldots
\Psi_{j_k}(\tau_{l_k})\sqrt{(\tau_{l_k})}
\Delta{\bf w}_{\tau_{l_k}}^{(i_k)}
\right)
+
$$

\vspace{2mm}
$$
+{\tilde R}_{T,t}^{p_1,\ldots,p_k}=
$$

\vspace{6mm}
$$
=\sum_{j_1=0}^{p_1}\ldots\sum_{j_k=0}^{p_k}
{\tilde C}_{j_k\ldots j_1}\times
$$

\vspace{2mm}
$$
\times
\left(
\prod_{l=1}^k {\tilde \zeta}_{j_l}^{(i_l)}-
\hbox{\vtop{\offinterlineskip\halign{
\hfil#\hfil\cr
{\rm l.i.m.}\cr
$\stackrel{}{{}_{N\to \infty}}$\cr
}} }\sum_{(l_1,\ldots,l_k)\in {\rm G}_k}
\Psi_{j_1}(\tau_{l_1})\sqrt{r(\tau_{l_1})}\Delta{\bf w}_{\tau_{l_1}}^{(i_1)}
\ldots
\Psi_{j_k}(\tau_{l_k})\sqrt{(\tau_{l_k})}
\Delta{\bf w}_{\tau_{l_k}}^{(i_k)}
\right)+
$$

\vspace{2mm}
$$
+{\tilde R}_{T,t}^{p_1,\ldots,p_k},
$$

\vspace{7mm}
\noindent
where

$$
{\tilde R}_{T,t}^{p_1,\ldots,p_k}
=\sum_{(t_1,\ldots,t_k)}
\int\limits_{t}^{T}
\ldots
\int\limits_{t}^{t_2}
\left(K(t_1,\ldots,t_k)\prod_{l=1}^k\sqrt{r(t_l)}-\right.
$$

\vspace{2mm}
$$
\left.
-\sum_{j_1=0}^{p_1}\ldots
\sum_{j_k=0}^{p_k}
{\tilde C}_{j_k\ldots j_1}
\prod_{l=1}^k\left(\Psi_{j_l}(t_l)\sqrt{r(t_l)}\right)\right)
d{\bf w}_{t_1}^{(i_1)}
\ldots
d{\bf w}_{t_k}^{(i_k)},
$$

\vspace{5mm}
\noindent
where permutations $(t_1,\ldots,t_k)$ when summing are performed only 
in the values $d{\bf w}_{t_1}^{(i_1)}
\ldots $
$d{\bf w}_{t_k}^{(i_k)}$. At the same time the indexes near 
upper limits of integration in the iterated stochastic integrals 
are changed correspondently and if $t_r$ swapped with $t_q$ in the  
permutation $(t_1,\ldots,t_k)$, then $i_r$ swapped $i_q$ in the 
permutation $(i_1,\ldots,i_k)$.

Let us evaluate the remainder
${\tilde R}_{T,t}^{p_1,\ldots,p_k}$ of the series.

According to Lemma 2 in \cite{arxiv-1} or Lemma 1.2 in 
\cite{2018a} (also see \cite{2018aa}-\cite{12aa-afterxxx}), we have

$$
{\sf M}\left\{\left({\tilde R}_{T,t}^{p_1,\ldots,p_k}\right)^2\right\}
\le 
C_k
\sum_{(t_1,\ldots,t_k)}
\int\limits_{t}^{T}
\ldots
\int\limits_{t}^{t_2}
\left(K(t_1,\ldots,t_k)\prod_{l=1}^k\sqrt{r(t_l)}\right.-
$$

\vspace{2mm}
\begin{equation}
\label{obana1eee}
\left.-
\sum_{j_1=0}^{p_1}\ldots
\sum_{j_k=0}^{p_k}
{\tilde C}_{j_k\ldots j_1}
\prod_{l=1}^k\left(\Psi_{j_l}(t_l)\sqrt{r(t_l)}\right)\right)^2
dt_1
\ldots
dt_k=
\end{equation}

\vspace{2mm}
$$
=C_k\int\limits_{[t,T]^k}
\left(K(t_1,\ldots,t_k)-
\sum_{j_1=0}^{p_1}\ldots
\sum_{j_k=0}^{p_k}
{\tilde C}_{j_k\ldots j_1}
\prod_{l=1}^k\Psi_{j_l}(t_l)\right)^2
\times
$$

\vspace{2mm}
$$
\times
\left(\prod_{l=1}^k r(t_l)\right)
dt_1 \ldots
dt_k\to 0
$$

\vspace{5mm}
\noindent
if $p_1,\ldots,p_k\to\infty,$ where constant $C_k$ 
depends only
on the multiplicity $k$ of the iterated Ito stochastic integral
(\ref{ito-rr}). 
Theorem 5 is proved.

Let us formulate the following theorem (the
version of Theorem 3 in \cite{arxiv-2}).

\vspace{2mm}

{\bf Theorem 6} \cite{2018a}-\cite{12aa-afterxxx}. 
{\it Suppose that
every $\psi_l(\tau)$ $(l=1,\ldots, k)$ is a continuous non-random function on 
$[t, T]$ and
$\{\Psi_j(x)\sqrt{r(x)}\}_{j=0}^{\infty}$ $(r(x)\ge 0)$
is a complete orthonormal 
system of functions in the space $L_2([t,T]),$ each function 
$\Psi_j(x)\sqrt{r(x)}$
of which 
for finite $j$ satisfies the condition 
$(\star)$ {\rm (}see Sect.~{\rm 1)}.
Then the estimate

\vspace{-2mm}
$$
{\sf M}\left\{\left(
{\tilde J}[\psi^{(k)}]_{T,t}-{\tilde J}[\psi^{(k)}]_{T,t}^{p_1,\ldots,p_k}
\right)^2\right\}
\le
$$

\begin{equation}
\label{z1-ura}
\le k!\left(\int\limits_{[t,T]^k}
K^2(t_1,\ldots,t_k)\left(\prod_{l=1}^k r(t_l)\right)
dt_1\ldots dt_k -\sum_{j_1=0}^{p_1}\ldots
\sum_{j_k=0}^{p_k}{\tilde C}^2_{j_k\ldots j_1}\right)
\end{equation}

\vspace{5mm}
\noindent
is valid for the following cases{\rm :}

\vspace{2mm}

{\rm 1.}\ $i_1,\ldots,i_k=1,\ldots,m$\ \ and\ \ $0<T-t<\infty,$

\vspace{1mm}

{\rm 2.}\ $i_1,\ldots,i_k=0, 1,\ldots,m,$\ \ $i_1^2+\ldots+i_k^2>0,$\ \
and\ \ $0<T-t<1,$

\vspace{2mm}

\noindent
where ${\tilde J}[\psi^{(k)}]_{T,t}$ is the 
stochastic integral {\rm (\ref{ito-rr}),}
${\tilde J}[\psi^{(k)}]_{T,t}^{p_1,\ldots,p_k}$ is the 
expression on the right-hand side of {\rm (\ref{tyyy-rr})} before
passing to the limit 
$\hbox{\vtop{\offinterlineskip\halign{
\hfil#\hfil\cr
{\rm l.i.m.}\cr
$\stackrel{}{{}_{p_1,\ldots,p_k\to \infty}}$\cr
}} };$ another 
notations are the same as in Theorem {\rm 5}.
}

\vspace{2mm}

Consider the following generalizations of Theorems 5, 6.

\vspace{2mm}

{\bf Theorem 7} \cite{2018a} (Sect.~1.13), \cite{arxiv-1} (Sect.~17).
{\it Let $\psi_1(x)\sqrt{r(x)},\ldots,
\psi_k(x)\sqrt{r(x)}\in L_2([t, T]),$ 
where $r(x)\ge 0.$
Furthermore$,$ let 
$\{\Psi_j(x)\sqrt{r(x)}\}_{j=0}^{\infty}$
is an arbitrary complete orthonormal
system of functions in the space $L_2([t,T]).$
Then$,$ for the iterated Ito stochastic integral

\vspace{-1mm}
\begin{equation}
\label{fifi1}
{\tilde J}[\psi^{(k)}]_{T,t}=\int\limits_t^T\psi_k(t_k)\sqrt{r(t_k)} 
\ldots \int\limits_t^{t_{2}}
\psi_1(t_1)\sqrt{r(t_1)} d{\bf w}_{t_1}^{(i_1)}\ldots
d{\bf w}_{t_k}^{(i_k)}
\end{equation}

\vspace{3mm}
\noindent
the following expansion 

$$
{\tilde J}[\psi^{(k)}]_{T,t}=
\hbox{\vtop{\offinterlineskip\halign{
\hfil#\hfil\cr
{\rm l.i.m.}\cr
$\stackrel{}{{}_{p_1,\ldots,p_k\to \infty}}$\cr
}} }
\sum\limits_{j_1=0}^{p_1}\ldots
\sum\limits_{j_k=0}^{p_k}
C_{j_k\ldots j_1}\Biggl(
\prod_{l=1}^k {\tilde \zeta}_{j_l}^{(i_l)}+\sum\limits_{r=1}^{[k/2]}
(-1)^r \times
\Biggr.
$$

\vspace{2mm}
\begin{equation}
\label{fifi2}
\times
\sum_{\stackrel{(\{\{g_1, g_2\}, \ldots, 
\{g_{2r-1}, g_{2r}\}\}, \{q_1, \ldots, q_{k-2r}\})}
{{}_{\{g_1, g_2, \ldots, 
g_{2r-1}, g_{2r}, q_1, \ldots, q_{k-2r}\}=\{1, 2, \ldots, k\}}}}
\prod\limits_{s=1}^r
{\bf 1}_{\{i_{g_{{}_{2s-1}}}=~i_{g_{{}_{2s}}}\ne 0\}}
\Biggl.{\bf 1}_{\{j_{g_{{}_{2s-1}}}=~j_{g_{{}_{2s}}}\}}
\prod_{l=1}^{k-2r} {\tilde \zeta}_{j_{q_l}}^{(i_{q_l})}\Biggr)
\end{equation}

\vspace{5mm}
\noindent
that converges in the mean-square
sense   
is valid, where 
$i_1,\ldots,i_k=0,1,\ldots,m,$ 

\vspace{-1mm}
$$
{\tilde \zeta}_{j}^{(i)}=
\int\limits_t^T \Psi_{j}(s)\sqrt{r(s)}d{\bf w}_s^{(i)}
$$

\vspace{3mm}
\noindent
are independent standard Gaussian random variables
for various
$i$ or $j$ {\rm(}in the case when $i\ne 0${\rm),}

$$
{\tilde C}_{j_k\ldots j_1}=\int\limits_{[t,T]^k}
K(t_1,\ldots,t_k)
\prod_{l=1}^{k}\biggl(\Psi_{j_l}(t_l)r(t_l)\biggr)dt_1\ldots dt_k
$$

\vspace{3mm}
\noindent
is the Fourier coefficient$,$
$K(t_1,\ldots,t_k)$ is defined by {\rm (\ref{ppp});}
another notations are the same as in Theorems {\rm 1, 2, 5.}
}

\vspace{2mm}

{\bf Theorem 8} \cite{2018a} (Sect.~1.13), \cite{arxiv-1} (Sect.~17).
{\it Let $\psi_1(x)\sqrt{r(x)},\ldots,
\psi_k(x)\sqrt{r(x)}\in L_2([t, T]),$ 
where $r(x)\ge 0.$
Furthermore$,$ let 
$\{\Psi_j(x)\sqrt{r(x)}\}_{j=0}^{\infty}$
is an arbitrary complete orthonormal 
system of functions in the space $L_2([t,T]).$
Then the following estimate

\vspace{1mm}
$$
{\sf M}\left\{\left(
{\tilde J}[\psi^{(k)}]_{T,t}-{\tilde J}[\psi^{(k)}]_{T,t}^{p_1,\ldots,p_k}
\right)^2\right\}
\le 
$$

\vspace{2mm}
$$
~ \le k!\left(~\int\limits_{[t,T]^k}
K^2(t_1,\ldots,t_k)\left(\prod_{l=1}^k r(t_l)\right)
dt_1\ldots dt_k -\sum_{j_1=0}^{p_1}\ldots
\sum_{j_k=0}^{p_k}{\tilde C}^2_{j_k\ldots j_1}\right)
$$

\vspace{5mm}
\noindent
is valid for the following cases{\rm :}

\vspace{2mm}

{\rm 1.}\ $i_1,\ldots,i_k=1,\ldots,m$\ \ and\ \ $0<T-t<\infty,$

\vspace{1mm}

{\rm 2.}\ $i_1,\ldots,i_k=0, 1,\ldots,m,$\ \ $i_1^2+\ldots+i_k^2>0,$\ \
and\ \ $0<T-t<1,$

\vspace{2mm}
\noindent
where ${\tilde J}[\psi^{(k)}]_{T,t}$ is the 
stochastic integral {\rm (\ref{fifi1}),}
${\tilde J}[\psi^{(k)}]_{T,t}^{p_1,\ldots,p_k}$ is the 
expression on the right-hand side of {\rm (\ref{fifi2})} before
passing to the limit 
$\hbox{\vtop{\offinterlineskip\halign{
\hfil#\hfil\cr
{\rm l.i.m.}\cr
$\stackrel{}{{}_{p_1,\ldots,p_k\to \infty}}$\cr
}} };$ another 
notations are the same as in Theorem {\rm 2, 7}.
}

\vspace{5mm}

\section{One Modification of Theorems 4 and 5}

\vspace{5mm}

Let us compare (\ref{obana1eee}) and (\ref{z2}). 
If we suppose that $r(x)\ge 0$ and 

\vspace{1mm}
$$
\frac{\rho(x)}{r(x)}\le C<\infty,
$$  

\vspace{2mm}
\noindent
where $\rho(x)$ as in (\ref{riss100}), then 

\vspace{-1mm}
$$
\int\limits_{[t,T]^k}
\left(K(t_1,\ldots,t_k)
-\sum_{j_1=0}^{p_1}\ldots
\sum_{j_k=0}^{p_k}
{\tilde C}_{j_k\ldots j_1}
\prod_{l=1}^k\Psi_{j_l}(t_l)\right)^2\times
$$

$$
\times
\rho(t_1)dt_1
\ldots
\rho(t_k)dt_k=
$$

\vspace{3mm}
$$
=
\int\limits_{[t,T]^k}
\left(K(t_1,\ldots,t_k)
-\sum_{j_1=0}^{p_1}\ldots
\sum_{j_k=0}^{p_k}
{\tilde C}_{j_k\ldots j_1}
\prod_{l=1}^k\Psi_{j_l}(t_l)\right)^2\times
$$

\vspace{1mm}
$$
\times
\frac{\rho(t_1)}{r(t_1)}r(t_1)dt_1
\ldots
\frac{\rho(t_k)}{r(t_k)}r(t_k)dt_k\le
$$

\vspace{3mm}
$$
\le
C_k'\int\limits_{[t,T]^k}
\left(K(t_1,\ldots,t_k)
-\sum_{j_1=0}^{p_1}\ldots
\sum_{j_k=0}^{p_k}
{\tilde C}_{j_k\ldots j_1}
\prod_{l=1}^k\Psi_{j_l}(t_l)\right)^2\times
$$

$$
\times
r(t_1)dt_1
\ldots
r(t_k)dt_k,
$$

\vspace{5mm}
\noindent
where $C_k'$ is a constant,
$\{\Psi_j(x)\}_{j=0}^{\infty}$ is a complete orthonormal 
with weight $r(x)\ge 0$ 
system of functions in the space $L_2([t, T]),$
and the Fourier coefficient 
${\tilde C}_{j_k\ldots j_1}$
has the form (\ref{koef}).

So, we obtain the following modification of Theorems 4 and 5.

\vspace{2mm}

{\bf Theorem 9}\ \cite{2018a}, \cite{arxiv-26b}.
{\it Suppose that the following 
conditions are fulfilled{\rm :}

{\rm 1}. Every $\psi_l(\tau)\ (l=1,\ldots,k)$ is a 
continuous non-random function at
the interval $[t, T]$.

{\rm 2}. $M_{\tau}^{(l,i_l)}\in Q_4(\rho,[t,T]),$ $G_n(\rho,[t,T])$ with
$n=2^{k+1},$  
$i_l=1,\ldots,m,$ $l=1,\ldots,k$ $(k\in\mathbb{N}).$

{\rm 3}. $\{\Psi_j(x)\}_{j=0}^{\infty}$ is a complete orthonormal 
with weight $r(\tau)\ge 0$ 
system of functions in the space $L_2([t,T]),$ each function of which 
for finite $j$ satisfies the condition $(\star)$ {\rm (}see Sect.~{\rm 1)}.
Moreover,

\vspace{-1mm}
$$
\frac{\rho(x)}{r(x)}\le C<\infty.
$$  
\vspace{1mm}

Then, for the iterated stochastic integral 
$J[\psi^{(k)}]_{T,t}^M$ with respect to martingales 
defined by {\rm (\ref{mart})}
the following expansion

\vspace{1mm}
$$
J[\psi^{(k)}]_{T,t}^M=
\hbox{\vtop{\offinterlineskip\halign{
\hfil#\hfil\cr
{\rm l.i.m.}\cr
$\stackrel{}{{}_{p_1,\ldots,p_k\to \infty}}$\cr
}} }\sum_{j_1=0}^{p_1}\ldots\sum_{j_k=0}^{p_k}
{\tilde C}_{j_k\ldots j_1}\Biggl(
\prod_{l=1}^k\xi_{j_l}^{(l,i_l)}-
\Biggr.
$$

\vspace{3mm}
$$
-\Biggl.
\hbox{\vtop{\offinterlineskip\halign{
\hfil#\hfil\cr
{\rm l.i.m.}\cr
$\stackrel{}{{}_{N\to \infty}}$\cr
}} }\sum_{(l_1,\ldots,l_k)\in {\rm G}_k}
\Psi_{j_{1}}(\tau_{l_1})
\Delta{M}_{\tau_{l_1}}^{(1,i_1)}\ldots
\Psi_{j_{k}}(\tau_{l_k})
\Delta{M}_{\tau_{l_k}}^{(k,i_k)}\Biggr)
$$

\vspace{6mm}
\noindent
converging in the mean-square sense is valid, 
where $i_1,\ldots,i_k=1,\ldots,m,$
$\left\{\tau_{j}\right\}_{j=0}^{N}$ is a partition of
the interval $[t, T]$ which satisfies the condition {\rm (\ref{w11}),}
$\Delta{M}_{\tau_{j}}^{(r,i)}=
M_{\tau_{j+1}}^{(r,i)}-M_{\tau_{j}}^{(r,i)}$
$(i=1,\ldots,m,$ $r=1,\ldots,k),$

\vspace{1mm}
$$
{\rm G}_k={\rm H}_k\backslash{\rm L}_k,\ \ 
{\rm H}_k=\{(l_1,\ldots,l_k):\ l_1,\ldots,l_k=0,\ 1,\ldots,N-1\},
$$

$$
{\rm L}_k=\{(l_1,\ldots,l_k):\ l_1,\ldots,l_k=0,\ 1,\ldots,N-1;\
l_g\ne l_r\ (g\ne r);\ g, r=1,\ldots,k\},
$$

\vspace{4mm}
\noindent
${\rm l.i.m.}$ is a limit in the mean-square sense$,$
$$
\xi_{j}^{(l,i_l)}=
\int\limits_t^T \Psi_{j}(s) d{M}_s^{(l,i_l)}
$$

\vspace{2mm}
\noindent
are independent for various 
$i_l=1,\ldots,m$ $(l=1,\ldots,k)$
and uncorrelated for various $j$
$\left(\hbox{if}\ i_l\ne 0,\right.$ 
$\left.\rho(x)\equiv r(x)\right)$ random variables,

$$
{\tilde C}_{j_k\ldots j_1}=\int\limits_{[t,T]^k}
K(t_1,\ldots,t_k)
\prod_{l=1}^{k}\biggl(\Psi_{j_l}(t_l)r(t_l)\biggr)dt_1\ldots dt_k
$$

\vspace{3mm}
\noindent
is the Fourier coefficient,

$$
K(t_1,\ldots,t_k)=
\begin{cases}
\psi_1(t_1)\ldots \psi_k(t_k)\ &\hbox{for}\ \ t_1<\ldots<t_k\\
~\\
~\\
0\ &\hbox{otherwise}
\end{cases},\ \ \ \ t_1,\ldots,t_k\in[t, T],\ \ \ \ k\ge 2
$$

\vspace{5mm}
\noindent
and 
$K(t_1)\equiv\psi_1(t_1)$ for $t_1\in[t, T].$  
}

\vspace{5mm}

{\bf Remark 2.}\ {\it Note that if $\rho(\tau), r(\tau)\equiv 1$ in 
Theorem {\rm 9},
then we obtain the variant of Theorem {\rm 1}.}

\vspace{5mm}

\section{Example on Application of Theorem 9 for the System
of Bessel Functions}

\vspace{5mm}

Let us consider the following boundary-value problem

$$
\left(p(x)\Phi'(x)\right)'+q(x)\Phi(x)=-\lambda r(x)\Phi(x),
$$

$$
\alpha\Phi(a)+\beta\Phi'(a)=0,\ \ \ 
\gamma\Phi(b)+\delta\Phi'(b)=0,
$$

\vspace{4mm}
\noindent
where the functions $p(x)$, $q(x)$, $r(x)$ satisfy 
the well known conditions and 
$\alpha,$ $\beta,$ $\gamma,$ $\delta,$ $\lambda$ are real numbers.

It is well known (Steklov V.A.) that 
the eigenfunctions
$\Phi_0(x),$ $\Phi_1(x),$ $\ldots $ of this boundary-value problem 
form a complete
orthonormal with weight $r(x)$ system of functions in the space 
$L_2([a, b]).$ 
It means that 
the Fourier series of the function $\sqrt{r(x)}f(x)\in L_2([a, b])$
with respect to the system of functions 

\vspace{-2mm}
$$
\sqrt{r(x)}\Phi_0(x),\ \ \ \sqrt{r(x)}\Phi_1(x),\ \ \ \ldots 
$$

\vspace{4mm}
\noindent
converges in the mean-square sense to the function $\sqrt{r(x)}f(x)$
at the interval $[a, b]$. Moreover,
the Fourier coefficients are defined by the formula

\vspace{-2mm}
\begin{equation}
\label{www.67}
C_j=\int\limits_a^b r(x)f(x)\Phi_j(x)dx.
\end{equation}

\vspace{3mm}

It is known that when solving the problem on oscillations of 
a circular membrane (general case), a boundary-value problem arises 
for the following Euler--Bessel equation

\vspace{-2mm}
\begin{equation}
\label{www.45}
r^2 R''(r)+rR'(r)+\left(\lambda^2 r^2-n^2\right)R(r)=0\ \ \ (\lambda\in\mathbb{R},\ \ \ 
n\in\mathbb{N}).
\end{equation}

\vspace{5mm}

The eigenfunctions of this problem, taking into account 
specific boundary 
conditions, are the following functions

\vspace{-2mm}
\begin{equation}
\label{www.55}
J_n\biggl(\mu_j\frac{r}{L}\biggr),
\end{equation}

\vspace{3mm}
\noindent
where $\tau\in[0, L]$ and $\mu_j$ $(j=0, 1, 2,\ldots )$ 
are positive roots of the Bessel 
function $J_n(\mu)$ ($n=0, 1, 2,\ldots$)
numbered in ascending order.

The problem on radial oscillations of a circular membrane leads to
the boundary-value problem for 
the
equation (\ref{www.45}) 
for $n=0$, 
the eigenfunctions of which are the functions (\ref{www.55}) when 
$n=0$.

Let us analyze the system of functions

\begin{equation}
\label{dmitri4}
\Psi_j(\tau)=\frac{\sqrt{2}}{T J_{n+1}(\mu_j)} J_n\left(
\frac{\mu_j}{T}\tau\right),\ \ \ j=0, 1, 2,\ldots,
\end{equation}

\vspace{4mm}
\noindent
where 

$$
J_n(x)=\sum\limits_{m=0}^{\infty}(-1)^m\left(\frac{x}{2}\right)^{n+2m}
\frac{1}{\Gamma(m+1)\Gamma(m+n+1)}
$$

\vspace{4mm}
\noindent
is the 
Bessel function of the first 
kind
and

$$
\Gamma(z)=\int\limits_0^{\infty}e^{-x} x^{z-1} dx
$$

\vspace{4mm}
\noindent
is the gamma-function, $\mu_j$ 
are positive  
roots of the function $J_n(x)$ numbered in ascending order, and
$n$ is a natural number or zero.

Due to the well known properties of 
the Bessel functions, the system $\left\{\Psi_j(\tau)\right\}_{j=0}^{\infty}$
is a complete orthonormal system of continuous functions 
with weight $\tau$ in the space $L_2([0, T])$.

Let us use the system of functions (\ref{dmitri4}) in 
Theorem 9.

Consider the following iterated stochastic integral with respect to
martingales

\vspace{-1mm}
$$
\int\limits_0^T\int\limits_0^s dM_{\tau}^{(1)}dM_s^{(2)},
$$

\vspace{3mm}
\noindent
where 
$$
M_s^{(i)}=\int\limits_0^s\sqrt{\tau}d{\bf f}_{\tau}^{(i)}\ \ \ (i=1, 2),
$$

\vspace{3mm}
\noindent
${\bf f}_{\tau}^{(i)}$ $(i=1, 2)$
are independent standard Wiener 
processes, 
$M_s^{(i)}$ $(i=1, 2)$ are martingales (here
$\rho(\tau)\equiv\tau$), $0\le s\le T.$
In addition, $M_s^{(i)}$ has a Gaussian distribution. 

It is obvious
that the conditions of Theorem 9 are fulfilled for $k=2$.
Using Theorem 9, we obtain

$$
\int\limits_0^T\int\limits_0^s dM_{\tau}^{(1)}dM_s^{(2)}=
\hbox{\vtop{\offinterlineskip\halign{
\hfil#\hfil\cr
{\rm l.i.m.}\cr
$\stackrel{}{{}_{p_1,p_2\to \infty}}$\cr
}} }\sum_{j_1=0}^{p_1}\sum_{j_2=0}^{p_2}
{\tilde C}_{j_2j_1}
\zeta_{j_1}^{(1)}\zeta_{j_2}^{(2)},
$$

\vspace{4mm}
\noindent
where 
$$
\zeta_j^{(i)}=\int\limits_0^T\Psi_j(\tau)dM_{\tau}^{(i)}
$$

\vspace{3mm}
\noindent
are independent standard  Gaussian random variables
for various $i$ or $j$ $(i=1, 2,$\ \ $j=0, 1, 2,\ldots),$

\vspace{1mm}
$$
{\sf M}\left\{\zeta_{j_1}^{(1)}\zeta_{j_2}^{(2)}\right\}=0,
$$

\vspace{2mm}
$$
{\tilde C}_{j_2 j_1}=
\int\limits_0^T s\Psi_{j_2}(s)
\int\limits_0^s \tau\Psi_{j_1}(\tau)d\tau ds.
$$

\vspace{4mm}

It is obvious that we can get this result using the another approach: 
we can use Theorems 1, 2 for the iterated Ito stochastic integral

\vspace{-1mm}
$$
\int\limits_0^T \sqrt{s}\int\limits_0^s \sqrt{\tau}d{\bf f}_{\tau}^{(1)}
d{\bf f}_s^{(2)},
$$

\vspace{3mm}
\noindent
and as a system of functions $\{\phi_j(s)\}_{j=0}^{\infty}$
in Theorems 1, 2 we can take

$$
\phi_j(s)=
\frac{\sqrt{2s}}{T J_{n+1}(\mu_j)} J_n\left(
\frac{\mu_j}{T}s\right),\ \ \ j=0, 1, 2,\ldots.
$$

\vspace{4mm}

As a result, we obtain

\vspace{-1mm}
$$
\int\limits_0^T \sqrt{s}\int\limits_0^s \sqrt{\tau}d{\bf f}_{\tau}^{(1)}
d{\bf f}_s^{(2)}=
\hbox{\vtop{\offinterlineskip\halign{
\hfil#\hfil\cr
{\rm l.i.m.}\cr
$\stackrel{}{{}_{p_1,p_2\to \infty}}$\cr
}} }\sum_{j_1=0}^{p_1}\sum_{j_2=0}^{p_2}
{C}_{j_2j_1}
\zeta_{j_1}^{(1)}\zeta_{j_2}^{(2)},
$$

\vspace{3mm}
\noindent
where 
$$
\zeta_j^{(i)}=\int\limits_0^T\phi_j(\tau)d{\bf f}_{\tau}^{(i)}
$$

\vspace{3mm}
\noindent
are independent standard Gaussian random variables
for various $i$ or $j$ $(i=1, 2,$\ \
$j=0, 1, 2,\ldots ),$

\vspace{1mm}
$$
{\sf M}\left\{\zeta_{j_1}^{(1)}\zeta_{j_2}^{(2)}\right\}=0,\ \ \ 
C_{j_2 j_1}=
\int\limits_0^T {\sqrt s}\phi_{j_2}(s)
\int\limits_0^s {\sqrt \tau}\phi_{j_1}(\tau)d\tau ds
$$

\vspace{3mm}
\noindent
is the Fourier coefficient. Obviously that 
$C_{j_2 j_1}={\tilde C}_{j_2 j_1}.$

Easy calculation demonstrates that

\vspace{1mm}
$$
\tilde\phi_j(s)=
\frac{\sqrt{2(s-t)}}{(T-t) J_{n+1}(\mu_j)} J_n\left(
\frac{\mu_j}{T-t}(s-t)\right),\ \ \ j=0, 1, 2,\ldots
$$ 

\vspace{4mm}
\noindent
is a complete orthonormal system of functions in the space 
$L_2([t, T]).$

Then, using Theorems 1, 2, we obtain

$$
\int\limits_t^T \sqrt{s-t}\int\limits_t^s \sqrt{\tau-t}d{\bf f}_{\tau}^{(1)}
d{\bf f}_s^{(2)}=
\hbox{\vtop{\offinterlineskip\halign{
\hfil#\hfil\cr
{\rm l.i.m.}\cr
$\stackrel{}{{}_{p_1,p_2\to \infty}}$\cr
}} }\sum_{j_1=0}^{p_1}\sum_{j_2=0}^{p_2}
{C}_{j_2j_1}
\tilde\zeta_{j_1}^{(1)}\tilde\zeta_{j_2}^{(2)},
$$

\vspace{4mm}
\noindent
where 
$$
\tilde\zeta_j^{(i)}=\int\limits_t^T\tilde\phi_j(\tau)d{\bf f}_{\tau}^{(i)}
$$

\vspace{3mm}
\noindent
are independent standard Gaussian random variables
for various $i$ or $j$ $(i=1, 2,$\ \
$j=0, 1, 2,\ldots ),$

\vspace{2mm}
$$
{\sf M}\left\{\tilde\zeta_{j_1}^{(1)}\tilde\zeta_{j_2}^{(2)}\right\}=0,\ \ \ 
C_{j_2 j_1}=
\int\limits_t^T \sqrt{s-t}\tilde\phi_{j_2}(s)
\int\limits_t^s \sqrt{\tau-t}\tilde\phi_{j_1}(\tau)d\tau ds.
$$


\vspace{12mm}


\begin{thebibliography}{199}



\vspace{10mm}



\bibitem{1}
Gihman I.I., Skorohod A.V. Stochastic Differential Equations and 
its Applications.
Kiev, Naukova Dumka, 1982, 612 pp. 



\bibitem{KlPl2}
Kloeden P.E., Platen E. Numerical Solution of Stochastic
Differential Equations. 
Berlin, Springer, 1992, 632 pp.



\bibitem{Mi2}
Milstein G.N. Numerical Integration of Stochastic Differential 
Equations. Sverdlovsk, Ural University Press, 1988, 225 pp. 


\bibitem{Mi3}
Milstein G.N., Tretyakov M.V. Stochastic Numerics
for Mathematical Physics. 
Berlin, Springer, 2004, 616 pp.



\bibitem{KPS}
Kloeden P.E., Platen E., Schurz H. Numerical Solution
of SDE Through Computer Experiments. Berlin, Springer, 1994, 292 pp.



\bibitem{PW1}
Platen E., Wagner W. On a Taylor formula for a class of Ito
processes. Probab. Math. Statist. 3 (1982), 37-51.


\bibitem{KlPl1}
Kloeden P.E., Platen E. The Stratonovich and Ito-Taylor 
Expansions. Math. Nachr. 151 (1991), 33-50.



\bibitem{kk5}
Kulchitskiy O.Yu., Kuznetsov D.F. The unified Taylor-Ito expansion.
Journal of Mathematical Sciences (N.~Y.). 99, 2  (2000), 1130-1140.	
DOI: https://doi.org/10.1007/BF02673635


\bibitem{kk6}
Kuznetsov D.F. New representations of the Taylor-Stratonovich expansions.
Journal of Mathematical Sciences (N.~Y.). 118, 6 (2003), 5586-5596.	
DOI: http://doi.org/10.1023/A:1026138522239



\bibitem{2006}
Kuznetsov D.F. Numerical Integration of Stochastic Differential Equations. 2.
[In Russian]. Polytechnical University Publishing House, 
Saint-Petersburg, 2006, 764 pp.
DOI: http://doi.org/10.18720/SPBPU/2/s17-227\
Available at:\ http://www.sde-kuznetsov.spb.ru/06.pdf\
(ISBN 5-7422-1191-0)




\bibitem{2011-2}
Kuznetsov D.F. Strong Approximation of Multiple Ito and 
Stratonovich Stochastic Integrals: Multiple Fourier Series Approach.
2nd Edition. [In English]. 
Polytechnical University Publishing House, Saint-Petersburg, 
2011, 284 pp. DOI: http://doi.org/10.18720/SPBPU/2/s17-233\
Available at:\\
http://www.sde-kuznetsov.spb.ru/11a.pdf\
(ISBN 978-5-7422-3162-2)




\bibitem{2017}
Kuznetsov D.F. Multiple Ito and Stratonovich Stochastic Integrals: 
Fourier-Legendre and Trigonometric Expansions, Approximations, Formulas.
[In English].
Electronic Journal "Differential Equations and Control Processes"
ISSN 1817-2172 (online),
1 (2017), A.1--A.385.\\
DOI: http://doi.org/10.18720/SPBPU/2/z17-3\\ 
Available at:\
http://diffjournal.spbu.ru/EN/numbers/2017.1/article.2.1.html 




\bibitem{2017-1}
Kuznetsov D.F. Stochastic Differential Equations: Theory and Practice of 
Numerical Solution. With Programs on MATLAB, 5th Edition. [In Russian].
Electronic Journal "Differential Equations and Control Processes"
ISSN 1817-2172 (online), 2 (2017), A.1-A.1000.
DOI: http://doi.org/10.18720/SPBPU/2/z17-4\
Available at:\\
http://diffjournal.spbu.ru/EN/numbers/2017.2/article.2.1.html



\bibitem{2018}
Kuznetsov D.F. Stochastic Differential Equations: Theory and Practice of 
Numerical Solution. With MATLAB Programs, 6th Edition. [In Russian].
Electronic Journal "Differential Equations and Control Processes"
ISSN 1817-2172 (online), 4 (2018), A.1-A.1073.\
Available at:\\
http://diffjournal.spbu.ru/EN/numbers/2018.4/article.2.1.html


                    
\bibitem{2018a}
Kuznetsov D.F.
Strong Approximation of Iterated Ito and Stratonovich Stochastic 
Integrals Based on Generalized Multiple Fourier Series. 
Application to Numerical Solution of Ito SDEs and Semilinear SPDEs.
[In English].
arXiv:2003.14184 [math.PR], 2023, 992 pp.


\bibitem{2018aa}
Kuznetsov D.F.
Strong Approximation of Iterated Ito and Stratonovich Stochastic 
Integrals Based on Generalized Multiple Fourier Series. 
Application to Numerical Solution of Ito SDEs and Semilinear SPDEs.
Electronic Journal "Differential Equations and Control Processes"
ISSN 1817-2172 (online), 4 (2020), A.1-A.606.\
Available at:\
http://diffjournal.spbu.ru/EN/numbers/2020.4/article.1.8.html


\bibitem{2018aaa}
Kuznetsov D.F.
Mean-Square Approximation of Iterated It\^{o} and Stratonovich Stochastic 
Integrals Based on Generalized Multiple Fourier Series. 
Application to Numerical Integration of It\^{o} SDEs and Semilinear SPDEs.
[In English].
Electronic Journal "Differential Equations and Control Processes"
ISSN 1817-2172 (online),
4 (2021), A.1-A.788.\ Available at:\
http://diffjournal.spbu.ru/EN/numbers/2021.4/article.1.9.html


\bibitem{12aa-afterxxx}
Kuznetsov, D.F.
Mean-Square Approximation of Iterated It\^{o} and Stratonovich Stochastic 
Integrals Based on Generalized Multiple Fourier Series. 
Application to Numerical Integration of It\^{o} SDEs and Semilinear SPDEs
(Third Edition).
[In English].
Electronic Journal "Differential Equations and Control Processes"
ISSN 1817-2172 (online),
1 (2023), A.1-A.947.\ Available at:\
http://diffjournal.spbu.ru/EN/numbers/2023.1/article.1.10.html




\bibitem{2007-1}
Kuznetsov D.F. Stochastic Differential Equations: Theory and Practice 
of Numerical Solution. With MatLab Programs, 1st Edition. [In Russian]. 
Polytechnical University Publishing House, Saint-Petersburg, 2007, 778 pp.
DOI: http://doi.org/10.18720/SPBPU/2/s17-228\
Available at:\ http://www.sde-kuznetsov.spb.ru/07b.pdf\
(ISBN 5-7422-1394-8)




\bibitem{2007-2}
Kuznetsov D.F. Stochastic Differential Equations: Theory and Practice 
of Numerical Solution. With MatLab Programs, 2nd Edition. [In Russian]. 
Polytechnical 
University Publishing House, Saint-Petersburg, 2007, XXXII+770 pp.
DOI: http://doi.org/10.18720/SPBPU/2/s17-229\
Available at:\\
http://www.sde-kuznetsov.spb.ru/07a.pdf\
(ISBN 5-7422-1439-1)




\bibitem{2009}
Kuznetsov D.F. Stochastic Differential Equations: Theory and Practice 
of Numerical Solution. With MatLab Programs, 3rd Edition. [In Russian]. 
Polytechnical 
University Publishing House, Saint-Petersburg, 2009, XXXIV+768 pp.
DOI: http://doi.org/10.18720/SPBPU/2/s17-230\
Available at:\\
http://www.sde-kuznetsov.spb.ru/09.pdf\
(ISBN 978-5-7422-2132-6)





\bibitem{2010-1}
Kuznetsov D.F. Stochastic Differential Equations: Theory and Practice 
of Numerical Solution. With MatLab Programs. 4th Edition. [In Russian].
Polytechnical University Publishing House, Saint-Petersburg, 2010, 
XXX+786 pp. DOI: http://doi.org/10.18720/SPBPU/2/s17-231\
Available at:\
http://www.sde-kuznetsov.spb.ru/10.pdf\
(ISBN 978-5-7422-2448-8)




\bibitem{2010-2}
Kuznetsov D.F. Multiple Stochastic Ito and Stratonovich Integrals 
and Multiple Fourier Series.
[In Russian].
Electronic Journal "Differential Equations and Control Processes"
ISSN 1817-2172 (online),
3 (2010), A.1-A.257. DOI: http://doi.org/10.18720/SPBPU/2/z17-7\
Available at:\\
http://diffjournal.spbu.ru/EN/numbers/2010.3/article.2.1.html





\bibitem{2011-1}
Kuznetsov D.F. Strong Approximation of Multiple Ito and 
Stratonovich Stochastic Integrals: Multiple Fourier Series Approach.
1st Edition. [In English]. 
Polytechnical University Publishing House, Saint-Petersburg, 
2011, 250 pp. DOI: http://doi.org/10.18720/SPBPU/2/s17-232\
Available at:\\
http://www.sde-kuznetsov.spb.ru/11b.pdf\
(ISBN 978-5-7422-2988-9)






\bibitem{2013}
Kuznetsov D.F. Multiple Ito and Stratonovich Stochastic 
Integrals: Approximations, Properties, Formulas. [In English].
Polytechnical University Publishing House, Saint-Petersburg,
2013, 382 pp.\\
DOI: http://doi.org/10.18720/SPBPU/2/s17-234\\
Available at:\
http://www.sde-kuznetsov.spb.ru/13.pdf\
(ISBN 978-5-7422-3973-4)






\bibitem{17a}
Kuznetsov D.F. Development and application of the Fourier 
method for the numerical solution of Ito stochastic differential 
equations. [In English]. Computational Mathematics and 
Mathematical Physics, 58, 7 (2018), 1058-1070.
DOI: http://doi.org/10.1134/S0965542518070096



\bibitem{18a}
Kuznetsov D.F. On numerical modeling of the multidimensional 
dynamic systems under random perturbations with the 1.5 and 2.0 
orders of strong convergence [In English]. Automation and Remote Control, 
79, 7 (2018), 1240-1254.
DOI: http://doi.org/10.1134/S0005117918070056



\bibitem{arxiv-14}
Kuznetsov D.F.
To numerical modeling with strong orders 1.0, 1.5, and 2.0 of 
convergence for multidimensional dynamical systems with random disturbances.
[In English].
arXiv:1802.00888 [math.PR]. 2018, 29 pp. 



\bibitem{28a}
Kuznetsov D.F. On numerical modeling of the multidimentional dynamic 
systems under random perturbations with the 2.5 order of strong 
convergence. [In English]. Automation and Remote Control, 
80, 5 (2019), 867-881. DOI: http://doi.org/10.1134/S0005117919050060



\bibitem{29a}
Kuznetsov D.F. Comparative analysis of the efficiency of application 
of Legendre polynomials and trigonometric functions to the numerical 
integration of It\^{o} stochastic differential equations.
[In English]. Computational Mathematics and Mathematical Physics, 
59, 8 (2019),  1236-1250.\\
DOI: http://doi.org/10.1134/S0965542519080116



\bibitem{30a}
Kuznetsov D.F. Expansion of multiple Stratonovich stochastic integrals 
of second multiplicity based on double Fourier-Legendre series 
summarized by Prinsheim method [In Russian]. 
Electronic Journal "Differential Equations and Control Processes"
ISSN 1817-2172 (online), 1 (2018), 1-34.\
Available at:\\
http://diffjournal.spbu.ru/EN/numbers/2018.1/article.1.1.html
 


\bibitem{31a}
Kuznetsov D.F. Application of the method of approximation of iterated 
Ito stochastic integrals based on generalized multiple Fourier 
series to the high-order strong numerical methods for non-commutative 
semilinear stochastic partial differential equations. 
[In English].
arXiv:1905.03724 [math.GM], 2019, 41 pp. 



\bibitem{200a}
Kuznetsov D.F. Application of the method of approximation 
of iterated stochastic Ito integrals based on generalized multiple 
Fourier series to the high-order strong numeri\-cal methods 
for non-commutative semilinear stochastic partial 
differential equations. [In English]. 
Electronic Journal "Differential Equations and Control Processes"
ISSN 1817-2172 (online), 3 (2019), 18-62.
Available at:\\
http://diffjournal.spbu.ru/EN/numbers/2019.3/article.1.2.html


\bibitem{200aa}
Kuznetsov D.F. Application of multiple Fourier-Legendre series 
to implementation of strong exponential Milstein and 
Wagner-Platen methods for non-commutative semilinear stochastic 
partial differential equations. [In English].
arXiv:1912.02612 [math.PR], 2019, 32 pp.


\bibitem{200aaa}
Kuznetsov D.F. Application of multiple Fourier-Legendre series 
to strong exponential Milstein and 
Wagner-Platen methods for non-commutative semilinear stochastic 
partial differential equations. [In English].
Electronic Journal "Differential Equations and Control Processes"
ISSN 1817-2172 (online),
3 (2020), 129-162.
Available at:\
http://diffjournal.spbu.ru/EN/numbers/2020.3/article.1.6.html


\bibitem{new-2023a}
Kuznetsov D.F. A new proof of the expansion of iterated Ito stochastic 
integrals with respect to the components of a multidimensional Wiener 
process based on generalized multiple Fourier series and 
Hermite polynomials. [In English]. arXiv:2307.11006 [math.PR]. 2023, 56 pp. 



\bibitem{300a}
Kuznetsov D.F. A method of expansion and approximation of repeated 
stochastic Stratonovich integrals based on multiple Fourier series 
on full orthonormal systems. [In Russian].
Electronic Journal "Differential Equations and Control Processes"
ISSN 1817-2172 (online),
1 (1997), 18-77.
Available at:\\
http://diffjournal.spbu.ru/EN/numbers/1997.1/article.1.2.html



\bibitem{400a}
Kuznetsov D.F. Problems of the Numerical Analysis of Ito Stochastic 
Differential Equations. 
[In Russian].
Electronic Journal "Differential Equations and Control Processes"
ISSN 1817-2172 (online),
1 (1998), 66-367.
Available at:\\
http://diffjournal.spbu.ru/EN/numbers/1998.1/article.1.3.html\
Hard Cover Edition: 1998, SPbGTU Publishing House, 
204 pp. (ISBN 5-7422-0045-5)


\bibitem{500a}
Kuznetsov D.F. Mean square approximation of solutions 
of stochastic differential 
equations using Legendres polynomials. [In English]. Journal of 
Automation and 
Information Sciences (Begell House), 2000, 32 (Issue 12), 69-86.
DOI: http://doi.org/10.1615/JAutomatInfScien.v32.i12.80 


\bibitem{600a}
Kuznetsov D.F. New representations of explicit one-step numerical 
methods for jump-diffusion stochastic differential equations. 
[In English]. Computational 
Mathematics and Mathematical Physics, 41, 6 (2001), 874-888.
Available at:\ 
http://www.sde-kuznetsov.spb.ru/01b.pdf



\bibitem{301a}
Kuznetsov D.F. Comparative analysis of the efficiency of application 
of Legendre polynomials and trigonometric functions to the numerical 
integration of Ito stochastic differential equations.
[In English].
arXiv:1901.02345 [math.GM], 2019, 40 pp. 



\bibitem{271a}
Kuznetsov D.F. Expansion of iterated Stratonovich stochastic integrals
based on generalized multiple Fourier series. 
[In English]. Ufa Mathematical Journal, 
11, 4 (2019), 49-77. 
DOI: http://doi.org/10.13108/2019-11-4-49\\
Available at:\
http://matem.anrb.ru/en/article?art\_id=604.



\bibitem{arxiv-1}
Kuznetsov D.F. 
Expansion of iterated Ito stochastic integrals of arbitrary 
multiplicity based on generalized multiple Fourier series
converging in the mean. [In English].
arXiv:1712.09746 [math.PR]. 2023, 143 pp. 



\bibitem{arxiv-2}
Kuznetsov D.F.
Exact calculation of the mean-square error in the method of 
approximation of iterated Ito stochastic integrals based on 
generalized multiple Fourier series. [In English].
arXiv:1801.01079 [math.PR]. 2018, 70 pp. 



\bibitem{arxiv-3}
Kuznetsov D.F. Mean-square approximation of iterated Ito and 
Stratonovich stochastic integrals of multiplicities 1 to 6 
from the Taylor-Ito and 
Taylor-Stratonovich expansions using Legendre polynomials.
[In English].
arXiv:1801.00231 [math.PR]. 2017, 106 pp. 





\bibitem{arxiv-4}
Kuznetsov D.F.
The hypotheses on expansions of iterated Stratonovich stochastic 
integrals of arbitrary multiplicity and their partial proof. 
[In English].
arXiv:1801.03195 [math.PR]. 
2023, 159 pp. 




\bibitem{arxiv-5}
Kuznetsov D.F. 
Expansions of iterated Stratonovich stochastic integrals 
based on generalized multiple Fourier series:
multiplicities 1 to 6 and beyond. 
[In English].
arXiv:1712.09516 [math.PR]. 2023, 221 pp.





\bibitem{arxiv-7}
Kuznetsov D.F. 
Expansion of iterated Stratonovich stochastic integrals of 
multiplicity 3 based on generalized multiple Fourier series
converging in the mean: general case of series summation.
[In English].
arXiv:1801.01564 [math.PR]. 2018, 66 pp. 




\bibitem{arxiv-8}
Kuznetsov D.F. 
Expansion of iterated Stratonovich stochastic integrals of 
multiplicity 2 based on double Fourier-Legendre series
summarized by Pringsheim method. [In English].
arXiv:1801.01962 [math.PR]. 2018, 49 pp. 






\bibitem{arxiv-12}
Kuznetsov D.F. 
Development and 
application of the Fourier method to the mean-square approximation 
of iterated Ito and Stratonovich stochastic integrals. 
[In English].
arXiv:1712.08991 [math.PR]. 2022, 57 pp. 



\bibitem{arxiv-23}
Kuznetsov D.F.
Expansion of iterated Stratonovich stochastic integrals 
of arbitrary multiplicity
based on generalized iterated Fourier series converging pointwise. 
[In English].
arXiv:1801.00784 [math.PR]. 2018, 77 pp. 



\bibitem{arxiv-24}
Kuznetsov D.F. 
Strong numerical methods of orders 2.0, 2.5, and 3.0 for 
Ito stochastic differential equations based on the unified 
stochastic Taylor expansions and multiple Fourier-Legendre series.
[In English].
arXiv:1807.02190 [math.PR]. 2018, 44 pp.




\bibitem{arxiv-26b}
Kuznetsov D.F. Expansion of iterated stochastic integrals with respect to 
martingale Poisson measures and with respect to martingales based on 
generalized multiple Fourier series. [In English].
arXiv:1801.06501 [math.PR].
2018, 40 pp.



\bibitem{arxiv-9}
Kuznetsov D.F. Expansion of iterated Stratonovich stochastic 
integrals of fifth and sixth multiplicity based on generalized multiple Fourier 
series. [In English]. arXiv:1802.00643 [math.PR]. 2023, 149 pp. 


\bibitem{OK}
Kuznetsov D.F. The proof of convergence with probability 1 
in the method of expansion 
of iterated Ito stochastic integrals based on generalized multiple 
Fourier series.
Electronic Journal "Differential Equations and Control Processes"
ISSN 1817-2172 (online),
2 (2020), 89-117.\\ Available at:\ 
http://diffjournal.spbu.ru/RU/numbers/2020.2/article.1.6.html



\bibitem{new-new-1}
Kuznetsov D.F., Kuznetsov M.D. Mean-square approximation of iterated 
stochastic integrals from strong exponential Milstein and Wagner--Platen 
methods for non-commutative semilinear SPDEs based on multiple 
Fourier--Legendre series. Recent Developments in Stochastic Methods and 
Applications. ICSM-5 2020. 
Springer Proceedings in Mathematics \& Statistics, vol 371, Eds. 
Shiryaev, A.N., Samouylov, K.E., Kozyrev, D.V.
Springer, Cham, 2021, pp. 17-32.\ DOI: http://doi.org/10.1007/978-3-030-83266-7\_2


\bibitem{new-art-1-xxy}
Kuznetsov D.F. A new approach to the series expansion of iterated 
Stratonovich stochastic integrals of arbitrary multiplicity with respect 
to components of the multidimensional Wiener process. [In English].
Electronic Journal "Differential Equations and Control Processes"
ISSN 1817-2172 (online), 2 (2022), 83-186. 
Available at:\\
http://diffjournal.spbu.ru/EN/numbers/2022.2/article.1.6.html



\bibitem{new-art-1xxys}
Kuznetsov D.F. A new approach to the series expansion of iterated 
Stratonovich stochastic integrals of arbitrary multiplicity with 
respect to components of the multidimensional Wiener process. II.
[In English].
Electronic Journal "Differential Equations and Control Processes"
ISSN 1817-2172 (online), 4 (2022), 135-194. 
Available at:\
http://diffjournal.spbu.ru/EN/numbers/2022.4/article.1.9.html





\bibitem{allen}
Allen E. Approximation of triple 
stochastic integrals through region subdivision. Communicat. 
in Appl. Anal. Special Tribute Issue to Prof. V. Lakshmikantham.
17 (2013), 355-366.


\bibitem{rrr}
Averina T.A., Prigarin S.M.
Calculation of stochastic integrals of Wiener processes. 
Preprint 1048. Novosibirsk, 
Institute of Computational 
Mathematics and Mathematical Geophysics of
Siberian Branch of the Russian Academy of Sciences, 1995, 15 pp. 


\bibitem{rr}
Prigarin S.M., Belov S.M. On one application of the Wiener 
process decomposition into series. Preprint 1107. Novosibirsk, 
Siberian Branch of the Russian Academy of Sciences, 1998, 16 pp. 
[In Russian].



\bibitem{KPW}
Kloeden P.E., Platen E., Wright I.W. The approximation of multiple 
stochastic integrals. Stochastic Analysis and Applications.
10, 4 (1992), 431-441. 


\bibitem{Zapad-9}
Platen E., Bruti-Liberati N. Numerical Solution of Stochastic 
Differential Equations
with Jumps in Finance. Springer, Berlin-Heidelberg, 2010. 868 pp.



\bibitem{W-Z-1}
Wong E., Zakai M. On the convergence of ordinary integrals to 
stochastic integrals. Ann. Math. Stat.,
5, 36 (1965), 1560-1564.


\bibitem{W-Z-2}
Wong E., Zakai M. On the relation between ordinary and stochastic 
differential equations. Int. J. Eng. Sci., 3 (1965), 213-229.


\bibitem{Watanabe}
Ikeda N., Watanabe S. Stochastic
Differential Equations and Diffusion Processes.
2nd Edition. North-Holland Publishing Company,
Amsterdam, Oxford, New-York, 1989. 555 pp.


\bibitem{Rybakov1000}
Rybakov K.A. Orthogonal expansion of multiple It\^{o} stochastic integrals.
Electronic Journal "Differential Equations and Control Processes"
ISSN 1817-2172 (online),
3 (2021), 109-140. Available at:\\
http://diffjournal.spbu.ru/EN/numbers/2021.3/article.1.8.html





\bibitem{Sc}
Skhorohod A.V. Stochastc Processes with Independent Augments. 
[In Russian]. Moscow, Nauka Publ., 1964. 280 pp.


\end{thebibliography}
\end{document}